%% file: nm_ent_s0.tex
\documentclass[final]{siamltex}

\input avz_defs

\title{Optimal Subharmonic Entrainment}

\author{Anatoly Zlotnik$^*$ \and Jr-Shin Li\thanks{Department of Electrical and Systems Engineering, Washington University in St. Louis ({\tt azlotnik@ese.wustl.edu, jsli@ese.wustl.edu}).}}

\begin{document}

\maketitle

\begin{abstract}
For many natural and engineered systems, a central function or design goal is the synchronization of one or more rhythmic or oscillating processes to an external forcing signal, which may be periodic on a different time-scale from the actuated process.  Such subharmonic synchrony, which is dynamically established when $N$ control cycles occur for every $M$ cycles of a forced oscillator, is referred to as $N$:$M$ entrainment.  In many applications, entrainment must be established in an optimal manner, for example by minimizing control energy or the transient time to phase locking.  We present a theory for deriving inputs that establish subharmonic $N$:$M$ entrainment of general nonlinear oscillators, or of collections of rhythmic dynamical units, while optimizing such objectives.  Ordinary differential equation models of oscillating systems are reduced to phase variable representations, each of which consists of a natural frequency and phase response curve.  Formal averaging and the calculus of variations are then applied to such reduced models in order to derive optimal subharmonic entrainment waveforms.  The optimal entrainment of a canonical model for a spiking neuron is used to illustrate this approach, which is readily extended to arbitrary oscillating systems.
\end{abstract}


%\begin{AMS}
%15A15, 15A09, 15A23
%\end{AMS}

\pagestyle{myheadings}
\thispagestyle{plain}
%\markboth{TEX PRODUCTION AND V. A. U. THORS}{SIAM MACRO EXAMPLES}

\section{Introduction} \label{secintro}

The synchronization of interacting cyclical processes that evolve on different time-scales plays a fundamental role in many natural phenomena and engineered structures \cite{blekhman88,pikovsky01}.  The concept of synchronization is particularly significant to the study of biological systems \cite{strogatz01,granada09enz}, which exhibit endogenous oscillations with periods ranging from milliseconds, such as in spiking neurons \cite{izhikevich07}, to years, as in hibernation cycles \cite{mrosovsky80}.  All life on earth consists of rhythmic systems that are affected by other cyclical processes, such as the daily light and temperature changes that actuate circadian pacemakers \cite{doyle06}, and the complex environmental cycles that drive plant growth \cite{mcclung11}.  Biological systems may also interact in networks of responsive dynamical units, such as the neurons that constitute circadian oscillators \cite{gonze05} and central pattern generators \cite{coombes05} in the brain, or communicating insects \cite{hartbauer05}, for which even the simplest models can produce very complex dynamics \cite{strogatz00}.

The process of entrainment, which refers to the dynamic synchronization of an oscillating system to a periodic input, is significant in biology \cite{hanson78,ermentrout84,glass88,granada09enz}, with particular relevance in neuroscience \cite{demir97,berke04,sirota08}, and is also observed in reactive chemical systems \cite{kuramoto84,aronson86,lev89}.  The notion of entrainment is paramount for understanding rhythmic systems, as well as for controlling such systems in an optimal manner \cite{granada09fast,harada10}.  Optimal entrainment also has compelling applications in clinical medicine, such as protocols for coping with jet lag \cite{vosko10,sthilaire12}, clinical treatments for neurological disorders including epilepsy \cite{kiss08,good09}, Parkinson's disease \cite{hofmann11}, and tinnitus \cite{strauss05}, and optimization of cardiac pacemakers \cite{naqvi12}.  Techniques for controlling the entrainment process can also be used in the design of vibrating mechanical structures \cite{blekhman88,zalalutdinov03} and nanoscale electromechanical devices \cite{feng08,barnes11} that require frequency control or phase locking, and can enable transformational technologies such as neurocomputers \cite{hoppensteadt00} and chaos communication \cite{fischer00}.  Various phenomena such as noise-induced synchronization \cite{wacker11}, time-scales in synchronization and network dynamics \cite{chen09,chavez05}, and transient phenomena \cite{granada09fast,zlotnik13prl} have been examined.

Nonlinear oscillating systems are often studied by transforming the complex dynamic equations that describe their behavior into phase coordinates \cite{ermentrout96,strogatz00}, which can also be experimentally established for a physical system when the dynamics are unknown \cite{galan05,tokuda07}.  Such models have been studied extensively, with a particular focus on neural \cite{ermentrout96,hoppensteadt99} and electrochemical \cite{kuramoto84,kiss02,nakata09} systems.  The reduction of a system model from a complicated set of differential equations to a simple scalar phase coordinate representation is especially compelling from a control-theoretic perspective because it enables a corresponding reduction in the complexity of optimal control problems involving that system.  Optimal control of phase models has been investigated with various objectives, such as to alter the spiking of a single neuron using minimum energy inputs \cite{moehlis06} with constrained amplitude \cite{dasanayake11,dasanayake11cdc} and charge balancing \cite{danzl10,dasanayake11cb}, as well as to control a network of globally coupled neurons \cite{nabi11}.  Several studies have focused on optimal waveforms for entrainment using basic models \cite{schaus06,harada10}, and our recent work has resulted in optimal entrainment controls for general nonlinear oscillators that require no knowledge about the initial state or phase of the system \cite{zlotnik11}, and can account for uncertainty in oscillation frequency \cite{zlotnik12jne}.  These investigations have demonstrated that phase coordinate reduction provides a practical approach to the optimal control of complex oscillating systems.

Previous work on optimal control of the entrainment process has focused on the harmonic case, which corresponds to a one-to-one (1:1) relationship between the frequencies of the stimulus and oscillator.  Many physical processes, however, undergo subharmonic $N$:$M$ entrainment, which transpires when $N$ cycles of the stimulus occur for every $M$ cycles of the oscillator \cite{glass88}.  Originally examined in the context of loudspeaker dynamics \cite{cunningham51}, subharmonic synchronization can emerge among weakly coupled oscillators \cite{ermentrout81,guevara82}, and can be induced in forced or injection-locked oscillators to produce entrainment \cite{daryoush89,storti88}.  Subharmonic locking phenomena are of interest in a wide range of fields, and neuroscience in particular.  Applications exist in magnetoencephalography \cite{tass98}, the study of brain connectivity \cite{honey07}, dynamic neural regulation \cite{hunter03}, as well as in the clinical treatment of epilepsy  \cite{hofmann11,kiss08,kriellaars94}.  Subharmonic entrainment plays a central role in our understanding of human perception of beat and meter \cite{large96,clayton05,nozaradan11}, as well as sound in general, and an ability to affect this phenomenon will lead to innovative therapies for tinnitus \cite{strauss05,tass12}.  Indeed, the functional connectivity of the cerebral cortex may be shaped by mutual entrainment of bursting neurons across multiple time scales in a coevolutionary manner \cite{honey07,lajoie11,gutierrez11}.  Previous studies have found that subharmonic synchronization phenomena are ubiquitous in biological systems.    In fact, respiration and heartbeat in human beings is typically entrained at a 1:4 ratio \cite{schafer98}, and evidence exists that human sleep latency is entrained by the lunar cycle  \cite{cajochen13,foster08}, which is a 1:28 ratio.  Other investigations have focused on engineering subharmonic locking in electronic circuits \cite{maffezzoni10,takano07}, antenna systems \cite{zarroug95}, and voice coil audio systems \cite{bolanos05,noris12}.

In this paper, we develop a method for engineering weak, periodic signals that achieve subharmonic entrainment in nonlinear oscillating systems without the use of state feedback.  We apply the methods of phase model reduction, formal averaging, and the calculus of variations, which have been used in our previous studies on harmonic entrainment \cite{zlotnik11,zlotnik12jne,zlotnik13prl}.  In addition to yielding optimal waveforms for entrainment using weak forcing, this approach allows us to approximate the entrainment regions called Arnold tongues whereby the frequency-locking characteristics of the controlled system are visualized \cite{ermentrout81,schaus06}.  We have previously used such graphs to characterize the performance of optimal controls derived using the phase response curve (PRC) of an oscillator when it is applied for harmonic entrainment of the original oscillator in state space \cite{zlotnik12jne}.  This crucial validation is extended to the methodology that we apply here to the subharmonic case.

In the following section, we discuss the phase coordinate transformation for a nonlinear oscillator and various methods for computing the PRC.  In Section \ref{secent}, we describe how averaging theory is used to study the asymptotic behavior of an oscillating system under subharmonic rhythmic forcing, and in Section \ref{secminpow} we use the calculus of variations to derive the minimum-energy subharmonic entrainment control for a single oscillator with arbitrary PRC.  In Section \ref{secfast}, a similar approach is applied to derive the control of fixed energy that produces the fastest subharmonic entrainment of a single oscillator.  In Section \ref{secnmens}, we formulate a result on minimum-energy subharmonic entrainment of ensembles of structurally similar nonlinear oscillators, and in Section \ref{secmaxr} we study the dual objective of entraining the largest collection of such oscillators with a control of fixed energy.  This is followed by Section \ref{sechh}, where we examine the performance of minimum energy subharmonic entrainment waveforms computed for the Hodgkin-Huxley model \cite{hodgkin52}, as well as Section \ref{secfastsim}, where the convergence rate of the system subject to inputs for fast subharmonic entrainment is examined.  Finally, in Section \ref{secdisc} we discuss several details and implications of this paper.  Throughout the manuscript, important concepts are described graphically using illustrations as well as examples using the phase model of the Hodgkin-Huxley system, which is described in Appendix \ref{aphh}, as a canonical nonlinear oscillator.

%subharmonic synch \cite{daryoush89,storti88}
%magnetoencephalography \cite{tass98}
%beat and meter and music \cite{large95,large96,clayton05,edward94,nozaradan11}
%epilepsy treatment \cite{hofmann11,kiss08,kriellaars94}
%antenna \cite{zarroug95}
%noise synchronization \cite{wacker11}
%tinnitus \cite{strauss05,tass12}
%cardiac \cite{guevara82}
%time-scales in synchronization and network dynamics \cite{chen09,chavez05}
%voice coil audio \cite{noris12,bolanos05}
%why study it \cite{wongsarnpigoon10}
%heartbeat and respiration 1:4 \cite{schafer98}
%fast entrainment \cite{granada09fast,zlotnik13prl}
%arnold tongues \cite{schaus06}

\section{Phase models} \label{secpv}

The phase coordinate transformation is a model reduction technique that is widely used for studying oscillating systems characterized by complex nonlinear dynamics \cite{kuramoto84}, and can also be used for system identification when the dynamics are unknown \cite{kiss05}.  Consider a full state-space model of an oscillating system, described by a smooth ordinary differential equation system
\begin{eqnarray} \label{sys1}
\Dx=f(x,u), \quad x(0)=x_0, \quad t\in[0,\infty)
\end{eqnarray}
where $x(t)\in\bR^n$ is the state and $u(t)\in\bR$ is a control.  Furthermore, we require that (\ref{sys1}) has an attractive, non-constant limit cycle $\gamma(t)=\gamma(t+T)$, satisfying $\dot{\gamma}=f(\gamma,0)$, on the periodic orbit $\Gamma=\set{y\in\bR^n}{y=\gamma(t) \text{ for } 0\leq t< T}\subset\bR^n$.    In order to study the behavior of this system, we reduce it to a scalar equation
\begin{eqnarray} \label{sys2}
\dot{\psi}=\w+Z(\psi)u,
\end{eqnarray}
which is called a phase model, where $Z$ is the phase response curve (PRC) and $\psi(t)$ is the phase associated to the isochron on which $x(t)$ is located.  The isochron is the manifold in $\bR^n$ on which all points have asymptotic phase $\psi(t)$ \cite{brown04}.  It is standard practice to define $\psi(t)=0$ (mod $2\pi$) when the first variable in the state vector $x$ attains its maximum over the orbit $\Gamma$.  This is due to the significant role of mathematical neuroscience in the development of phase model theory.  For models of neural oscillators, the first state variable often denotes the membrane potential, which exhibits spiking or relaxation behavior, so that $\psi(t)=2\pi k$ for $k=1,2,\ldots$ occur concurrently with successive spikes.  The conditions for validity and accuracy of phase reduced models have been determined \cite{efimov10,efimov11}, and the reduction is accomplished through the well-studied process of phase coordinate transformation \cite{efimov09}, which is based on Floquet theory \cite{perko90,kelley04}.  The model is assumed valid for inputs $u(t)$ such that the solution $x(t,x_0,u)$ to (\ref{sys1}) remains within a neighborhood $U$ of $\Gamma$.
\begin{figure}[h]
\centerline { \includegraphics[width=1.04\linewidth]{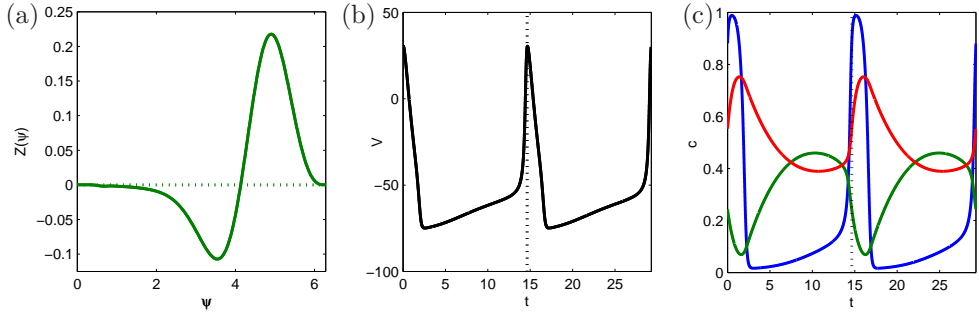}}
\begin{picture}(0,0)
\put(2,122){(a)} \put(129,122){(b)} \put(258,122){(c)}
\end{picture}
\vspace{-.5cm}
\caption{\footnotesize (a) Hodgkin-Huxley phase response curve (PRC). The natural period and frequency of oscillation are $T\approx 14.638$ ms and $\w\approx68.315$ rad/s, respectively.  (b) Voltage limit cycles. (c) Ion concentration limit cycles.  The motion along the periodic orbit is highly relaxational, as seen in the voltage ``spiking'' behavior.  The phase $\psi=0$ corresponds to ``spikes'', or maxima, of the membrane potential.} \vspace{.25cm}\label{s1f1_prc}
\end{figure}

To compute the PRC, the period $T=2\pi/\w$ and the limit cycle $\gamma(t)$ must be approximated to a high degree of accuracy.  This can be done using a method for determining the steady-state response of nonlinear oscillators \cite{aprille72} based on perturbation theory \cite{khalil02} and gradient optimization \cite{peressini00}.
The PRC can then be computed by integrating the adjoint of the linearization of (\ref{sys1}) \cite{ermentrout96}, or by using a more efficient and numerically stable spectral method developed more recently \cite{govaerts06}.  A software package called XPPAUT \cite{ermentrout02} is commonly used by researchers to compute the PRC.  We employ a technique derived from the method of Malkin \cite{malkin49} in order to compute PRCs, for which details are given in Appendix \ref{cprc}.  The PRC of the Hodgkin-Huxley system with nominal parameters, along with the limit cycle, is shown in Figure \ref{s1f1_prc}.

\section{Fundamental theory of subharmonic entrainment by weak forcing} \label{secent}

An essential objective in all entrainment applications is to force the frequency of an oscillator to a desired value.  While this can be accomplished using any sufficiently powerful rhythmic signal, it is often desirable to do so using an input that consumes minimum energy, or satisfies another optimization objective.  The harmonic ($1$:$1$) case constitutes the canonical entrainment problem, which was examined for arbitrary nonlinear oscillating systems in our previous work \cite{zlotnik11,zlotnik12jne}.  The theory of subharmonic ($N$:$M$) entrainment that is presented here is a nontrivial extension of those results.

Our goal is to entrain the system (\ref{sys2}) to a target frequency $\O$ using a periodic forcing control $u(t)$ of frequency $\Omega_f$, such that $M$ cycles of the oscillator occur for every $N$ cycles of the input.  When such $N$:$M$ entrainment occurs, then the target and forcing frequencies satisfy $M\Omega_f=N\Omega$, so that the control input has the form $u(t)=v(\frac{N}{M}\O t)$, where $v$ is $2\pi$-periodic.   From here on, it is assumed that $N$ and $M$ are coprime integers.  In addition, we adopt the weak forcing assumption, i.e., $v=\eP v_1$ where $v_1$ has unit energy and $\eP<<1$, so that given this control the state of the original system (\ref{sys1}) is guaranteed to remain in a neighborhood $U$ of $\Gamma$ in which the phase model (\ref{sys2}) remains valid \cite{efimov10}.  We then define a slow phase variable by $\phi(t)=\psi(t)-\O t$, and call the difference $\Delta\w=\w-\O$ between the natural and target frequencies the frequency detuning.  The dynamic equation for the slow phase is then
\begin{eqnarray} \label{sys3}
\dot{\phi}=\dot{\psi}-\O=\Delta\w+Z(\O t+\phi)v\bp{\frac{N}{M}\O t},
\end{eqnarray}
where $\dot{\phi}$ is called the phase drift.  In order to study the asymptotic behavior of (\ref{sys3}) it is necessary to eliminate the explicit dependence on time on the right hand side, which can be accomplished by using formal averaging \cite{kuramoto84}.  If $\cP$ is the set of $2\pi$-periodic functions on $\bR$, we can define an averaging operator $\bt{\cdot}:\cP\to\bR$ by
\begin{eqnarray} \label{ave}
\bt{x}=\frac{1}{2\pi}\int_0^{2\pi}x(\theta)\rd \theta.
\end{eqnarray}
In addition, let us define the forcing phase $\eta=\O_f t = \frac{N}{M}\O t$ and a change of variables $\theta=\eta/N$.  Then the weak ergodic theorem for measure-preserving dynamical systems on the torus \cite{kornfeld82} implies that for any periodic function $v$, the interaction function
\begin{align} \label{Lambdef}
\Lambda_v^\uNM(\phi)  & \teQ \bt{Z(M\theta+\phi)v(N\theta)} \notag \\& = \frac{1}{2\pi}\int_0^{2\pi} Z(M\theta+\phi)v(N\theta)\rd \theta \notag   \\   & =\lim_{T\to\infty}\frac{1}{T}\int_0^TZ(\O t+\phi)v\bp{\frac{N}{M}\O t}\rd t
\end{align}
exists as a continuous, $2\pi$-periodic function in $\cP$.   In addition, because both $Z$ and $v$ are $2\pi$-periodic, $\Lambda_v^{\uNM}$ can be expressed by integrating with respect to $\eta$ or to $\xi=M\theta=\O t$ to yield two equivalent expressions given by
\begin{align}
\Lambda_v^{\uNM}(\phi) &= \frac{1}{2\pi N}\int_0^{2\pi N} Z\bp{\frac{M}{N}\eta + \phi}v(\eta)\rd\eta  \notag %\label{lnm2}
\\ & = \frac{1}{2\pi N}\sum_{j=0}^{N-1}\int_0^{2\pi} Z\bp{\frac{M}{N}[2\pi j+\eta]+\phi}v(\eta)\rd \eta, \label{lnm2a} \\
\Lambda_v^{\uNM}(\phi) &= \frac{1}{2\pi M}\int_0^{2\pi M} Z\bp{\xi + \phi}v\bp{\frac{N}{M}\xi}\rd\xi \notag %\label{lnm3}
\\ & = \frac{1}{2\pi M}\sum_{\ell=0}^{M-1}\int_0^{2\pi} Z(\xi + \phi)v\bp{\frac{N}{M}[2\pi\ell+\xi]}\rd \xi. \label{lnm3a}
\end{align}
In particular, the expression \eqref{lnm2a} can be written as
\begin{align} \label{lnm2b}
\Lambda_v^{\uNM}(\phi) = \bt{Y^{\,\uNM}(\eta,\phi)v(\eta)}
\end{align}
where we define the function
\begin{align} \label{ynmdf}
Y^{\,\uNM}(\eta,\phi)= \frac{1}{N}\sum_{j=0}^{N-1}  Z\bp{\frac{M}{N}[2\pi j+\eta]+\phi}.
\end{align}
We henceforth write $Y^{\,\uNM}(\eta)\teQ Y^{\,\uNM}(\eta,0)$.  At this point, let us establish several important expressions that will be used throughout the following sections.  First, we define a function $Q$ as the $1$:$1$ interaction function of $Z$ with itself by
\begin{eqnarray} \label{qdef}
Q(\phi) \teQ \bt{Z(\theta+\phi)Z(\theta)}.
\end{eqnarray}
By defining an inner product $(\cdot,\cdot):\cP\times\cP\to\bR$ by $(f,g)=\bt{fg}$, the Cauchy-Schwarz inequality yields $|Q(\phi)|\leq\bt{Z^2}=Q(0)$, and the periodicity of $Z$ results in $Q(\phi)=\bt{Z(\theta+\phi)Z(\theta)} = \bt{Z(\theta)Z(\theta-\phi)} = Q(-\phi)$. We can then define
\begin{eqnarray} \label{vnmdef}
V^{\,\uNM}(\phi) \teQ \frac{1}{N}\sum_{j=0}^{N-1} Q\bp{\frac{M}{N}2\pi j+\phi},
\end{eqnarray}
which inherits the properties $|V^{\,\uNM}(\vphi)|\leq V^{\,\uNM}(0)$ for all $\vphi\in[0,2\pi)$ and $V^{\,\uNM}(-\vphi)=V^{\,\uNM}(\vphi)$ from the properties of $Q$.  We will subsequently write $V_0^{\,\uNM}\teQ V^{\,\uNM}(0)$ and  $V_*^{\,\uNM}=\min_{\phi\in[0,2\pi]}V^{\,\uNM}(\phi)$.  The expression \eqref{vnmdef} is important because using $v(\theta)=Y^{\,\uNM}(\theta,\psi)$ in \eqref{lnm3a} yields
\begin{align} \label{lnmv1}
\Lambda_{Y^{\,\uNM}(\theta,\psi)}^{\uNM}(\phi) & = \frac{1}{2\pi M}\sum_{\ell=0}^{M-1}\int_0^{2\pi} Z(\xi + \phi)Y^{\,\uNM}\bp{\frac{N}{M}[2\pi\ell+\xi],\psi}\rd \xi
 \notag \\
&=  \frac{1}{2\pi N}\sum_{j=0}^{N-1}\int_0^{2\pi} Z(\xi + \phi)Z\bp{\xi + \frac{M}{N}2\pi j+\psi}\rd \xi  \notag \\ &= \frac{1}{N}\sum_{j=0}^{N-1} Q\bp{\frac{M}{N}2\pi j+\phi-\psi} = V^{\,\uNM}(\phi-\psi).
\end{align}
In addition, using \eqref{lnm2b} we see that the energy of the function $Y^{\,\uNM}$ is given by
\begin{align} \label{lnmv2}
\bt{Y^{\,\uNM}Y^{\,\uNM}} & = \Lambda_{Y^{\,\uNM}}^{\uNM}(0) = V_0^{\,\uNM}.
\end{align}
The functions $Y^{\,\uNM}$, $Q$, and $V^{\,\uNM}$, as defined in \eqref{ynmdf}, \eqref{qdef}, and \eqref{vnmdef}, respectively, will appear repeatedly in the subsequent derivations of optimal subharmonic entrainment controls.

\begin{figure}[t]
\hspace{29pt} \centerline {\includegraphics[width=.7\linewidth]{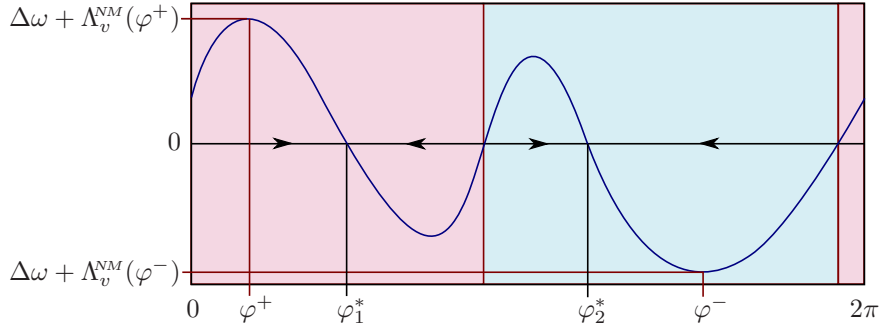}}
\begin{picture}(0,0)
\put(22,115){$\Delta\w+\Lambda_v^{\uNM}(\vphi^+)$} \put(22,19){$\Delta\w+\Lambda_v^{\uNM}(\vphi^-)$}
\put(109,5){$\vphi^+$} \put(281,5){$\vphi^-$} \put(146,5){$\vphi_1^*$} \put(237,5){$\vphi_2^*$}
\put(89,5){$0$} \put(340,5){$2\pi$} \put(82,68){$0$}
\end{picture}
\caption{\footnotesize Illustration of important properties of an interaction function $\Lambda_v^{\uNM}(\vphi)$.  The maximum and minimum values $\Lambda_v^{\uNM}(\vphi^+)$ and $\Lambda_v^{\uNM}(\vphi^-)$, which occur at the phases $\vphi^+$ and $\vphi^-$, respectively, determine the range of frequency detuning for which the oscillator can be entrained using weak forcing.  The roots of the equation $\Delta\w+\Lambda_v^{\uNM}(\vphi)=0$ determine the average phase shift, relative to $\O_f t$, at which the oscillation stabilizes from a given initial phase.  For initial phase in the pink (resp. blue) range, the asymptotic phase shift will be $\vphi_1^*$  (resp. $\vphi_2^*$).  The arrows indicate the evolution of the phase $\vphi$ in equation \ref{sys5}.} \vspace{.25cm} \label{s3f1_inter1}
\end{figure}

As in the case of 1:1 entrainment, the formal averaging theorem \cite{hoppensteadt97} permits us to approximate (\ref{sys3}) by the averaged system
\begin{eqnarray} \label{sys4}
\dot{\vphi}=\Delta \w + \Lambda_v^{\uNM}(\vphi)+\cO(\eP^2)
\end{eqnarray}
in the sense that there exists a change of variables $\vphi=\phi+\eP h(\vphi,\phi)$ that maps solutions of (\ref{sys3}) to those of (\ref{sys4}).  A detailed derivation for the $1$:$1$ case is given in Appendix B of \cite{zlotnik12jne}, and this can be easily extended to the $N$:$M$ case.  Therefore the weak forcing assumption $v=\eP v_1$ with $\eP<<1$ allows us to approximate the phase drift equation by
\begin{eqnarray} \label{sys5}
\dot{\vphi}=\Delta \w + \Lambda_v^{\uNM}(\vphi).
\end{eqnarray}
The averaged equation (\ref{sys5}) is autonomous, and approximately characterizes the asymptotic behavior of the system (\ref{sys2}) under periodic forcing.  Specifically, we say that the system is entrained by a control $u=v(\O_f t)$ when the phase drift equation (\ref{sys5}) satisfies $\dot{\vphi}=0$, which will occur as $t\to\infty$ if there exists a phase $\vphi_*$ that satisfies $\Delta \w + \Lambda_v^{\uNM}(\vphi^*) = 0$.  When both the control waveform $v$ and PRC $Z$ are non-zero, the function $\Lambda_v^{\uNM}(\vphi)$ is not identically zero, so when the system is entrained there exists at least one $\vphi_*\in[0,2\pi)$ that is an attractive fixed point of (\ref{sys5}).  The stable fixed points $\{\vphi_i^*\}$ of (\ref{sys5}), which are the roots of the equation $\Delta\w+\Lambda_v^{\uNM}(\vphi)=0$, determine the average phase shift, relative to $\O_f t$, at which the oscillation stabilizes from a given initial phase.  In addition, we define the phases $\vphi^+=\arg\max_\vphi\Lambda_v^{\uNM}(\vphi)$ and $\vphi^-=\arg\min_\vphi \Lambda_v^{\uNM}(\vphi)$ at which the interaction function achieves its maximum and minimum values, respectively.  In order for entrainment to occur, $-\Lambda_v^{\uNM}(\vphi^+)\leq \Delta\w \leq -\Lambda_v^{\uNM}(\vphi^-)$ must hold, so that at least one stable fixed point of $\Lambda_v^{\uNM}$ exists.  Thus the range of the interaction function determines which values of the frequency detuning $\Delta \w$ yield phase locking.  These properties are illustrated in Figure \ref{s3f1_inter1}.

\begin{figure}[t]
\hspace{10pt} \centerline {\includegraphics[width=.92\linewidth]{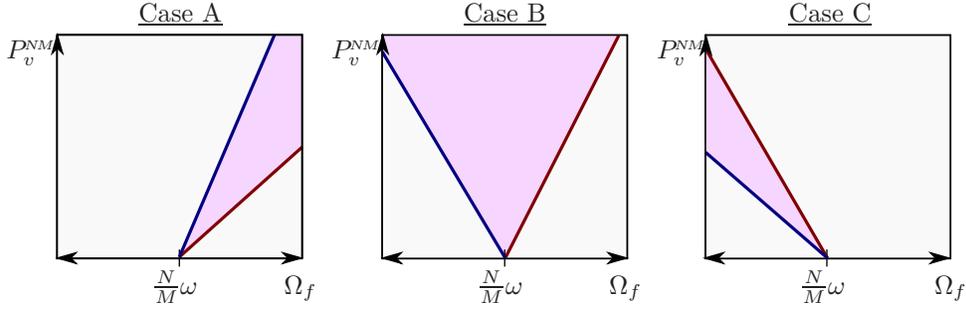}}
\begin{picture}(0,0)
\put(10,91){$P_v^{\,\uNM}$} \put(65,2){$\frac{N}{M}\w$} \put(115,2){$\O_f$} \put(60,105){\underline{Case A}}
\put(133,91){$P_v^{\,\uNM}$} \put(188,2){$\frac{N}{M}\w$} \put(238,2){$\O_f$} \put(183,105){\underline{Case B}}
\put(256,91){$P_v^{\,\uNM}$} \put(310,2){$\frac{N}{M}\w$} \put(361,2){$\O_f$} \put(306,105){\underline{Case C}}
\end{picture}
\caption{\footnotesize Illustration of Arnold tongues for the three cases listed in Table \ref{s3t1_arnold}.  The {\color{darkblue}left} boundary is shown in {\color{darkblue}blue}, and the {\color{darkred}right} boundary is shown in {\color{darkred}red}.} \vspace{.25cm} \vspace{.25cm} \label{s3f2_arnold}
\end{figure}

\renewcommand{\arraystretch}{1.2}
\begin{table}[t]
\centerline{
\begin{tabular}{|p{1.55cm}||c|c|c|c|}
\hline Case A: & \multicolumn{4}{|c|}{$0<\Lambda_{\tw{v}}^{\uNM}(\vphi^-)<\Lambda_{\tw{v}}^{\uNM}(\vphi^+)$}  \\ \hline Frequency: & \multicolumn{2}{|c|}{$\O_f>\frac{N}{M}\w$} & \multicolumn{2}{|c|}{$\O_f<\frac{N}{M}\w$} \\ \hline Boundary:  & {\color{darkblue}left/top} & {\color{darkred}right/bottom} & \multicolumn{2}{|c|}{N/A}  \\ \hline $P_v^{\,\uNM}(\O_f)$  & $-\Delta\w/\Lambda_{\tw{v}}^{\uNM}(\vphi^-)$ & $-\Delta\w/\Lambda_{\tw{v}}^{\uNM}(\vphi^+)$ &  \multicolumn{2}{|c|}{N/A} \\ \hline \hline
Case B: &  \multicolumn{4}{|c|}{$\Lambda_{\tw{v}}^{\uNM}(\vphi^-)<0<\Lambda_{\tw{v}}^{\uNM}(\vphi^+)$}  \\ \hline Frequency: & \multicolumn{2}{|c|}{$\O_f>\frac{N}{M}\w$} & \multicolumn{2}{|c|}{$\O_f<\frac{N}{M}\w$}  \\ \hline Boundary:  & \multicolumn{2}{|c|}{{\color{darkblue}left}} & \multicolumn{2}{|c|}{{\color{darkred}right}} \\ \hline $P_v^{\,\uNM}(\O_f)$  & \multicolumn{2}{|c|}{$-\Delta\w/\Lambda_{\tw{v}}^{\uNM}(\vphi^-)$} & \multicolumn{2}{|c|}{$-\Delta\w/\Lambda_{\tw{v}}^{\uNM}(\vphi^+)$} \\ \hline \hline
Case C: & \multicolumn{4}{|c|}{$\Lambda_{\tw{v}}^{\uNM}(\vphi^-)<\Lambda_{\tw{v}}^{\uNM}(\vphi^+)<0$} \\ \hline Frequency: & \multicolumn{2}{|c|}{$\O_f>\frac{N}{M}\w$} & \multicolumn{2}{|c|}{$\O_f<\frac{N}{M}\w$} \\ \hline Boundary: & \multicolumn{2}{|c|}{N/A}  & {\color{darkblue}left/bottom} & {\color{darkred}right/top} \\ \hline $P_v^{\,\uNM}(\O_f)$ &  \multicolumn{2}{|c|}{N/A} & $-\Delta\w/\Lambda_{\tw{v}}^{\uNM}(\vphi^-)$  & $-\Delta\w/\Lambda_{\tw{v}}^{\uNM}(\vphi^+)$ \\ \hline
\end{tabular}} \vspace{.25cm}
\caption{\footnotesize Arnold tongue boundary estimates for $N$:$M$ entrainment derived from \eqref{arnold_tg1}, where $\Delta \w \teQ \w-\frac{M}{N}\O_f$. } \vspace{.25cm} \label{s3t1_arnold}
\end{table}
\renewcommand{\arraystretch}{1}

Moreover, the interaction function can be used to estimate the values of the minimum root mean square (RMS) energy $P_v^{\,\uNM}(\O_f)=\sqrt{\bt{v^2}}$ that results in locking of an oscillator to a given frequency $\O_f$ at a subharmonic $N$:$M$ ratio using the waveform $v$.   This is accomplished by substituting the expression $v(\theta)=P_v^{\,\uNM}(\O_f)\tw{v}(\theta)$ and the relation $\O_f=\frac{N}{M}\O$ into the equation $\Delta\w+\Lambda_v^{\uNM}(\vphi)=0$ and simplifying to obtain \begin{align}
\w-\frac{M}{N}\O_f+P_v^{\,\uNM}(\O_f)\cdot\Lambda_{\tw{v}}^{\uNM}(\vphi)=0,
\label{arnold_tg1}
\end{align}
where $\tw{v}$ is a unit energy normalization of $v$.  This equation is then solved for $P_v^{\,\uNM}(\O_f)$ at $\vphi=\vphi^+$ and $\vphi=\vphi^-$ to produce linear estimates of boundaries for the regions of pairs $(\O_f,P_v^{\,\uNM})\in\bR^2$ that yield entrainment.   These regions are known as Arnold tongues, so named after mathematician who first described a similar phenomenon for recurrent maps on the circle (Section 12 of \cite{arnold61}). The RMS energy is used because the boundary of the entrainment region is approximately linear for weak forcing, and yields a clear visualization \cite{ermentrout81,schaus06}.  The Arnold tongue boundary estimates obtained using \eqref{arnold_tg1} can be classified into three different cases that depend on the signs of $\Lambda_{\tw{v}}^{\uNM}(\vphi^+)$ and $\Lambda_{\tw{v}}^{\uNM}(\vphi^-)$, which are listed in Table \ref{s3t1_arnold} and illustrated in Figure \ref{s3f2_arnold}.

Based on the theoretical foundation and fundamental notations presented in this section, we proceed to formulate and solve several design and optimization problems for subharmonic entrainment of one or more oscillating systems.  In the following section, we address the canonical problem of establishing subharmonic resonance of a single oscillator to a periodic input of minimum energy at a desired frequency.

\section{Minimum energy subharmonic entrainment of an oscillator} \label{secminpow}

In many applications described in Section \ref{secintro}, it is desirable to achieve entrainment of an oscillator by using a control of minimum energy.  This problem can be formulated as a variational optimization problem where the objective function to be minimized is the control energy $\bt{v^2}$, and the design constraint is $\w+\Lambda_v^{\uNM}(\vphi^+)\geq \O$ if $\O>\w$ and $\w+\Lambda_v^{\uNM}(\vphi^-)\leq\O$ if $\O<\w$.  This inequality is active when optimal entrainment occurs, and hence can be expressed as the equality constraint
\begin{eqnarray}
\Delta \w+\Lambda_v^{\uNM}(\vphi^+) &=& 0  \quad \text{if}\quad  \O>\w,\label{const1a}\\
\Delta \w+\Lambda_v^{\uNM}(\vphi^-) &=& 0 \quad \text{if}\quad \O<\w. \label{const1b}
\end{eqnarray}
We formulate the problem for $\O>\w$ to obtain the minimum energy control for frequency increase $v_+$  using the calculus of variations \cite{gelfand00}.  The derivation of the case where $\O<\w$ is similar, and results in the symmetric control $v_-$. The constraint (\ref{const1a}) can be adjoined to the cost  $\bt{v^2}$ using a multiplier $\lambda$, leading to the objective\\
\begin{align} \label{op1}
\cJ[v] & =  \bt{v^2} - \lambda(\Delta\w + \Lambda_v^{\uNM}(\vphi^+))  \\ & =
\bt{v^2} - \lambda\Delta\w \notag  -\frac{\lambda}{2\pi}\int_0^{2\pi} Y^{\,\uNM}(\eta,\vphi^+)v(\eta)\rd \eta\notag \\
& =  \frac{1}{2\pi} \int_0^{2\pi} \bp{ v(\eta)\bq{v(\eta) -\lambda Y^{\,\uNM}(\eta,\vphi^+)} - \lambda\Delta \w } \rd \eta,  \notag
\end{align}
where the expression \eqref{lnm2b} is substituted for $\Lambda_v^{\uNM}$.  Applying the Euler-Lagrange equation \cite{gelfand00}, we obtain the necessary condition for a candidate optimal solution
\begin{eqnarray} \label{sol0}
v_m(\eta)= \frac{\lambda}{2}Y^{\,\uNM}(\eta,\vphi^+).
\end{eqnarray}
Recalling \eqref{lnmv1}, we obtain
\begin{align} \label{lnm4}
\Lambda_{v_m}^{\uNM}(\vphi) & = \frac{\lambda}{2} \Lambda_{Y^{\,\uNM}(\theta,\vphi^+)}^{\uNM}(\phi) = \frac{\lambda}{2} V^{\,\uNM}(\vphi-\vphi^+).
\end{align}

\begin{figure}[t]
\centerline {\includegraphics[width=\linewidth]{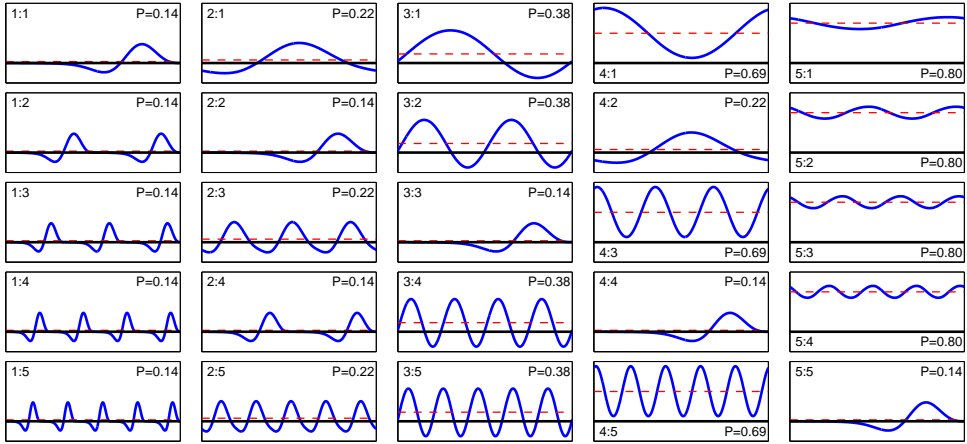}} \vspace{-.25cm}
\caption{\footnotesize Minimum energy subharmonic entrainment controls $v_+$ for increasing the frequency of the Hodgkin-Huxley neuron model by 3\%.  The controls for $N,\,M=1,\ldots,5$ are shown, and the domain and range in each plot are $[0,2\pi]$ and $[-.4,1.2]$, respectively.  The black line indicates the x-axis, and the red dashed line is the average value of the control.  The subharmonic ratio and RMS control power are indicated. Observe that if $N=1$, the control is simply repeated $M$ times to effectively produce 1:1 entrainment, which requires the lowest energy.  As $N$ grows large for a fixed $M$, the controls converge to a constant $\Upsilon$ given in equation \eqref{inf_const1}.} \vspace{.25cm} \label{s4f1_mpk}
\end{figure}

\begin{figure}[t]
\centerline {\includegraphics[width=\linewidth]{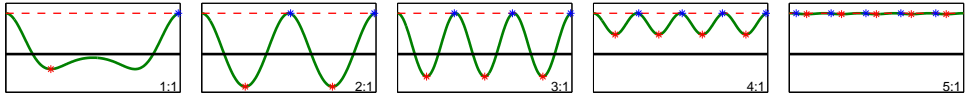}} \vspace{-.25cm}
\caption{\footnotesize Interaction functions $\Lambda_{v}^{\uNM}$ for minimum energy subharmonic entrainment controls $v_+$ for increasing the frequency of the Hodgkin-Huxley neuron model by 3\%, where the domain and range in each plot are $[0,2\pi]$ and $[-0.012,0.016]$, respectively.  Interaction functions for subharmonic entrainment depend only on $N$.  The black line indicates the x-axis, and the red dashed line is placed at $\Lambda_{v_+}^{\uNM}(\vphi^+)$, which is equal to the frequency detuning $\Delta\w=\w-\O$ by design.  As $N$ grows large, the interaction function converges to small variation about $\Delta\w$.  Red and blue stars mark the phases $\{\vphi^-_i\}$ and stable fixed points $\{\vphi_i^*\}$ of \eqref{sys5}, each of which occur $N$ times. } \vspace{.25cm} \label{s4f1_mpif}
\end{figure}

Therefore the constraint \eqref{const1a} results in
\begin{eqnarray} \label{lamnm1}
0 = \Delta \w+\Lambda_{v_m}^{\uNM}(\vphi^+) = \Delta \w + \frac{\lambda}{2}\, V_0^{\,\uNM},
\end{eqnarray}
so the multiplier is given by $\lambda=-2\Delta \w/V_0^{\,\uNM}$.  Consequently, we can express the minimum-energy control that entrains the system \eqref{sys2} to a target frequency $\O$ using a periodic forcing control $u(t)=v(\Omega_f t)$, where $v$ is $2\pi$-periodic and $\Omega_f=\frac{N}{M}\Omega$, as
\begin{align} \label{solnm1}
v_m(\eta)& = \left\{\begin{array}{ll} \dS v_+(\eta)=-\frac{\Delta \w}{V_0^{\,\uNM}} \cdot Y^{\,\uNM}(\eta,\vphi^+) &  \quad \text{if}\quad  \O>\w,\\\\\dS v_-(\eta)=-\frac{\Delta \w}{V_0^{\,\uNM}} \cdot Y^{\,\uNM}(\eta,\vphi^-)&  \quad \text{if}\quad  \O<\w, \end{array}\right.
\end{align}
where $\eta=\Omega_f t$ is the forcing phase. In practice, we may omit the phase ambiguity $\vphi^+$ or $\vphi^-$ in \eqref{solnm1} %and use $Y^{\,\uNM}(\eta)\teQ Y^{\,\uNM}(\eta,0)$
because entrainment is asymptotic.  The optimal waveforms for minimum-energy entrainment of the Hodgkin-Huxley system are shown in Figure \ref{s4f1_mpk} for values of $N,\,M=1,\ldots,5$, and the corresponding interaction functions are shown in Figure \ref{s4f1_mpif}.  Observe that in Figure \ref{s4f1_mpk}, the $1$:$M$ minimum-energy control will repeat the $1$:$1$ optimal waveform $M$ times during the control cycle.  As the ratio $N/M$ grows large, the controls converge to a constant given by
\begin{align}
\Upsilon=-\Delta\w\cdot\frac{\int_{0}^{2\pi}Z(\theta)\rd \theta}{\int_{0}^{2\pi}Q(\vphi)\rd \vphi},
\label{inf_const1}
\end{align}
as seen by substituting the limiting expressions as $N\to\infty$ for $Y^{\,\uNM}$ and  $V^{\,\uNM}$ from equations \eqref{ynmdf} and \eqref{vnmdef} into the solution \eqref{solnm1}.

Finally, by invoking \eqref{lnmv2}, we see that the minimum energy waveform \eqref{solnm1} for subharmonic entrainment of a single oscillator has energy given by
\begin{align} \label{minpvpow}
\bt{v_m^2}=\bp{\frac{\Delta\w}{V_0^{\,\uNM}}}^2\bt{Y^{\,\uNM}(\eta,\vphi^+)Y^{\,\uNM}(\eta,\vphi^+)}=\frac{(\Delta\w)^2}{V_0^{\,\uNM}}.
\end{align}

We have shown that the minimum energy periodic control $u(t)=v(\frac{N}{M}\O t)$ that achieves subharmonic entrainment of a single oscillator with natural frequency $\w$ to a target frequency $\O$ is a re-scaling of the function $Y^{\,\uNM}(\eta)$ given in \eqref{ynmdf}, where $\eta=\frac{N}{M}\O t$ is the forcing phase.  In the case that $N=M=1$, these results reduce to the solution in the harmonic (1:1) case, which is a re-scaling of the PRC $Z$ \cite{zlotnik11,zlotnik12jne}.  In addition, it is important to note that very similar optimal inputs for altering the frequency of oscillating neural systems in the phase model representation have been obtained using other methods in the harmonic (1:1) case \cite{moehlis06,dasanayake11}.  The fundamental conclusion is that when the desired change in the frequency of the oscillator is small and the applied input to the system is weak, the optimal control in the harmonic case is a re-scaling of the PRC.  Therefore we expect that applying alternative methods \cite{moehlis06,dasanayake11} to compute inputs for minimum energy subharmonic control of oscillators will also result in solutions similar to \eqref{solnm1}.

\section{Fast subharmonic entrainment of an oscillator} \label{secfast}

An alternative objective to minimizing control energy is to entrain a system to a desired frequency as quickly as possible using a control of given energy.  This problem has been examined for the harmonic ($1$:$1$) case for arbitrary nonlinear oscillating systems in our previous work \cite{zlotnik13prl}, and the solution for the subharmonic ($N$:$M$) case extends those results by applying the techniques derived in Section \ref{secent}.

\begin{figure}[t]
\centerline {\includegraphics[width=.7\linewidth]{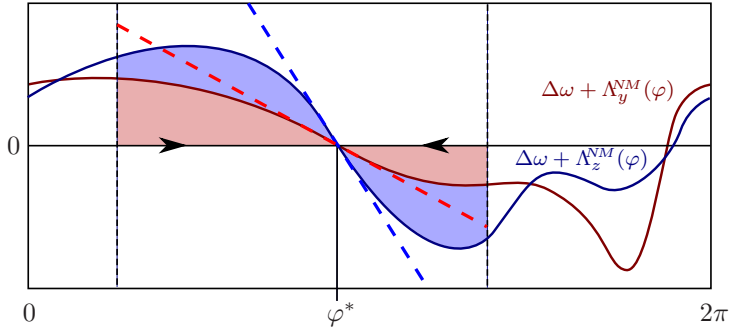}}
\begin{picture}(0,0)
\put(169,5){$\vphi^*$} \put(48,68){$0$} \put(240,64){{\color{darkblue}\footnotesize$\Delta\w+\Lambda_z^{\uNM}(\vphi)$}}
\put(250,90){{\color{darkred}\footnotesize$\Delta\w+\Lambda_y^{\uNM}(\vphi)$}}
\put(54,5){$0$} \put(310,5){$2\pi$}
\end{picture}
\caption{\footnotesize Illustration explaining the objective function for fast entrainment.  The averaged equation \eqref{sys5} is shown for two control waveforms {\color{darkred}$y$} and {\color{darkblue}$z$} that both result in the same attractive fixed phase $\vphi^*$.  Observe that ${\color{blue}-\ddx{\vphi}\Lambda_z(\vphi^*)}> {\color{red}-\ddx{\vphi}\Lambda_y(\vphi^*)}$, as indicated by the dashed lines.  As a result, $|\dot{\vphi}|$ is greater in the phase region between the dashed lines for the waveform $z$, as indicated by the shaded regions, so that the system converges to $\vphi^*$ faster when forced using $z$.  Therefore, we maximize the objective $\cJ[v]=-\ddx{\vphi}\Lambda_v(\vphi^*)$ for fast entrainment in problem \eqref{prob:1}.}  \vspace{.25cm} \label{s5f0_inter1}
\end{figure}

Our goal here is to entrain the system \eqref{sys2} to a target frequency $\O$ as quickly as possible by using a periodic control $v$ of fixed energy $P=\bt{v^2}$ and forcing frequency $\Omega_f$ that satisfies $M\Omega_f=N\Omega$.  Employing averaging theory as in Section \ref{secent} yields the phase drift equation \eqref{sys5}, where the interaction function would ideally be of a piecewise-constant form, so that the averaged slow phase $\vphi$ converges to a fixed point $\vphi^*$ at a uniform rate from any initial value.  However, a discontinuity at $\vphi\to\vphi^*$ would result in an unbounded control $v$, as explained in Lemma 2 of Appendix \ref{apnum}, which makes such a control infeasible in practice.  An alternative is to maximize $|\dot{\vphi}_*|$,  the rate of convergence of the averaged slow phase in the neighborhood of its attractive fixed point $\vphi^*$.  The calculus of variations can then be used to obtain a smooth optimal candidate solution that also performs well in practice.  When the system \eqref{sys5} is entrained by a control $v$, there exists an attractive fixed point $\vphi^*$ satisfying $\Lambda_v^{\uNM}(\vphi^*)+\Delta\w=0$ and $\ddx{\vphi}\Lambda_v^{\uNM}(\vphi^*)<0$, as seen in Figure \ref{s5f0_inter1}.   Observe that by inspecting \eqref{lnm2b}, one can write
\begin{align} \label{dlnmdf}
\ddx{\vphi}\Lambda_v^{\uNM}(\vphi^*)=\ddx{\vphi}\bt{Y^{\,\uNM}(\eta,\vphi^*)v(\eta)} = \bt{Y_\vphi^{\,\uNM}(\eta,\vphi^*)v(\eta)},
\end{align}
where $Y^{\,\uNM}$ is as defined in \eqref{ynmdf} and $Y_\vphi^{\,\uNM}(\eta,\vphi)$ is its derivative, given by
\begin{align} \label{dynmdf}
Y_\vphi^{\,\uNM}(\eta,\vphi)=\ddx{\vphi}Y^{\,\uNM}(\eta,\vphi) = & \frac{1}{N}\sum_{j=0}^{N-1}  Z'\bp{\frac{M}{N}[2\pi j+\eta]+\vphi},
\end{align}
with $Z'(\theta+\vphi)=\ddx{\vphi}Z(\theta+\vphi)$.  As was done for $Z$ in Section \ref{secent}, we can define the interaction function of $Z'$ with itself by
\begin{eqnarray} \label{rrdef}
K(\vphi) \teQ \bt{Z'(\theta+\vphi)Z'(\theta)},
\end{eqnarray}
which is maximized at $\vphi=0$ with the maximum value  $K(0)=\bt{Z'Z'}$.  The periodicity of $Z$ implies that  $|K(\vphi)|\leq\bt{Z'Z'}$ for all $\vphi\in[0,2\pi)$, and $K(\vphi)=K(-\vphi)$.   We then define
\begin{eqnarray} \label{wwnmdef}
S^{\,\uNM}(\vphi) \teQ \frac{1}{N}\sum_{j=0}^{N-1} K\bp{\frac{M}{N}2\pi j+\vphi},
\end{eqnarray}
which inherits the properties $|S^{\,\uNM}(\vphi)|\leq S^{\,\uNM}(0)$ for all $\vphi\in[0,2\pi)$ and $S^{\,\uNM}(-\vphi)=S^{\,\uNM}(\vphi)$ from the function $K$.  We will use the notation $S_0^{\,\uNM}\teQ S^{\,\uNM}(0)$ and  $S_*^{\,\uNM}=\min_{\phi\in[0,2\pi]}S^{\,\uNM}(\phi)$.   By substituting $Z'$ for $Z$ and $Y_\vphi^{\,\uNM}$ for $Y^{\,\uNM}$ in \eqref{lnmv1}, we can obtain
\begin{align} \label{dlnmv1}
\ddx{\phi}\Lambda_{Y_\vphi^{\,\uNM}(\theta,\psi)}^{\uNM}(\phi) & =  \frac{1}{N}\sum_{j=0}^{N-1} K\bp{\frac{M}{N}2\pi j+\phi-\psi} = S^{\,\uNM}(\phi-\psi).
\end{align}
In addition, by combining \eqref{dlnmdf} and \eqref{dlnmv1}, the energy of the function $Y_\vphi^{\,\uNM}$ is given by
\begin{align} \label{dlnmv2}
\bt{Y_\vphi^{\,\uNM}Y_\vphi^{\,\uNM}} & = \Lambda_{Y_\vphi^{\,\uNM}}^{\uNM}(0) = S_0^{\,\uNM}.
\end{align}
The functions $K$, and $S^{\,\uNM}$, as defined in \eqref{rrdef} and \eqref{wwnmdef}, respectively, will appear repeatedly in the following derivation of fast subharmonic entrainment controls.

In order to maximize the rate of entrainment in a neighborhood of $\vphi^*$ using a control of energy $P$, the value of $|\dot{\vphi}|$ should be maximized for values of $\vphi$ near $\vphi^*$, which occurs when $-\ddx{\vphi}\Lambda_v^{\uNM}(\vphi^*)$ is large, as illustrated in Figure \ref{s5f0_inter1}.  This results in the following optimal control problem formulation for fast subharmonic entrainment:
\begin{align}
	\label{prob:1}
	\max_{v\in\cP} \quad & \cJ[v] = -\ddx{\vphi}\Lambda_v^{\uNM}(\vphi^*)\\
	\label{eq:con1}
	\rM{s.t.} \quad & \bt{v^2} = P  \\
	\label{eq:con2}
	& \Lambda_v^{\uNM}(\vphi^*) + \Delta \w = 0.  %\\ & \Lambda_v'(\vphi^*) < 0 \nonumber
\end{align}
The constraints can be adjoined to the objective using multipliers $\lambda$ and $\mu$ to yield
\begin{align}
	\cJ[v] & =  -\ddx{\vphi}\Lambda_v^{\uNM}(\vphi^*) + \lambda(\bt{v^2} - P) + \mu(\Lambda_v^{\uNM}(\vphi^*)+\Delta\w) \nonumber\\
	& =  -\bt{Y_\vphi^{\,\uNM}(\eta,\vphi^*)v(\eta)} + \lambda(\bt{v^2}-P)  +  \mu\bt{ Y^{\,\uNM}(\eta,\vphi^*)v(\eta)} + \mu\Delta\w \nonumber\\
	&= \frac{1}{2\pi}\int_0^{2\pi} \bp{v(\eta)\bq{\mu Y^{\,\uNM}(\eta,\vphi^*) - Y_\vphi^{\,\uNM}(\eta,\vphi^*) + \lambda v(\eta)}   - \lambda P + \mu\Delta\w}\rd \eta.
\end{align}
The associated Euler-Lagrange equation is
\begin{align}
\mu Y^{\,\uNM}(\eta,\vphi^*) - Y_\vphi^{\,\uNM}(\eta,\vphi^*) + 2\lambda v(\eta) = 0,
\end{align}
and solving for $v$ yields the candidate solution
\begin{align} \label{sol:cand}
v_f(\eta)=\frac{1}{2\lambda}\bq{Y_\vphi^{\,\uNM}(\eta,\vphi^*) - \mu Y^{\,\uNM}(\eta,\vphi^*) }.
\end{align}
The multipliers $\lambda$ and $\mu$ can be found by substituting %the candidate solution
(\ref{sol:cand}) into the constraints \eqref{eq:con1} and \eqref{eq:con2}. %in problem (\ref{prob:1}).
This yields the equations
\begin{align} \label{const:1}
	\frac{1}{4\lambda^2}\bq{\bt{Y_\vphi^{\,\uNM}Y_\vphi^{\,\uNM}}-2\mu\bt{Y_\vphi^{\,\uNM} Y^{\,\uNM}} + \mu^2\bt{Y^{\,\uNM}Y^{\,\uNM}}} &= P, \\
	\label{const:2}
	\frac{1}{2\lambda}\bq{\bt{Y_\vphi^{\,\uNM} Y^{\,\uNM}} - \mu \bt{Y^{\,\uNM}Y^{\,\uNM}}} &= -\Delta \w,
\end{align}
where averaging is done with respect to the variable $\eta$.  Because $Z$ is $2\pi$-periodic, then $Z'$ is as well, as are $Y^{\,\uNM}$ and $Y_\vphi^{\,\uNM}$ in both arguments.  Thus one can show, e.g., using Fourier series, that $\bt{Z'Z}=0$, and $\bt{Y_\vphi^{\,\uNM} Y^{\,\uNM}}=0$ also, so that (\ref{const:2}) easily yields
\begin{align} \label{mufasteq}
\mu = \frac{2\Delta \w\lambda}{V_0^{\,\uNM}},
\end{align}
where $V_0^{\,\uNM}=\bt{Y^{\,\uNM}Y^{\,\uNM}}$ can be seen from \eqref{lnmv2}.  Substituting this result into (\ref{const:1}) leads to a quadratic equation for $\lambda$ given by
\begin{align} \label{lfasteq}
	\frac{1}{4\lambda^2}S_0^{\,\uNM} + \frac{(\Delta\w)^2}{V_0^{\,\uNM} } &= P.
\end{align}

\begin{figure}[t]
\centerline {\includegraphics[width=\linewidth]{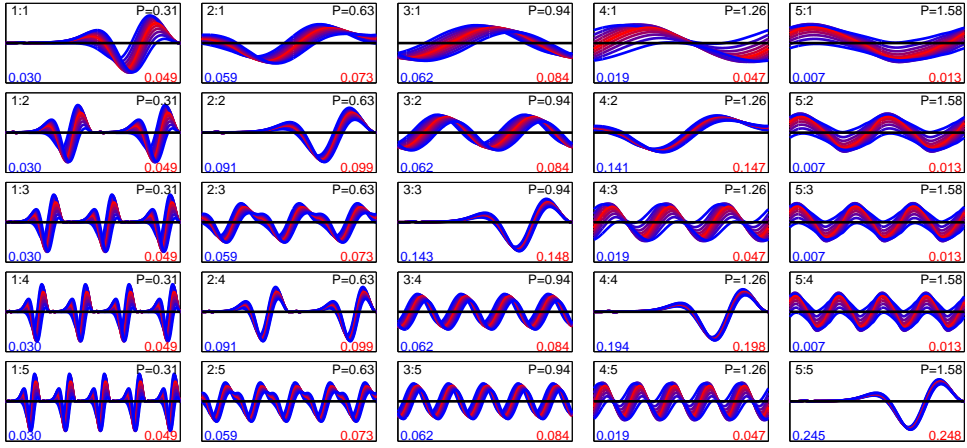}} \vspace{-.25cm}
\caption{\footnotesize Fast subharmonic entrainment controls $v_f$ for the Hodgkin-Huxley neuron model with frequency detuning of {\color{red}$0\%$} to {\color{blue}$\pm5\%$}.  The controls for $N,\,M=1,\ldots,5$ are shown rescaled to unit power, with the domain and range in each plot at $[0,2\pi]$ and $[-3.7,3.7]$, respectively.  The black line indicates the $x$-axis.  The entrainment ratio and RMS control energy are indicated; higher power is necessary to achieve a given detuning as the ratio $N$:$M$ increases. The entrainment rate, which is characterized by the slope of $\Lambda_{v_f}^{\uNM}(\vphi^*)$, is noted {\color{blue} in blue at bottom left for maximum detuning at $\pm5\%$}, and {\color{red}in red at bottom right for zero detuning}.} \vspace{.25cm} \label{s5f1_fk}
\end{figure}

\begin{figure}[t]
\centerline {\includegraphics[width=\linewidth]{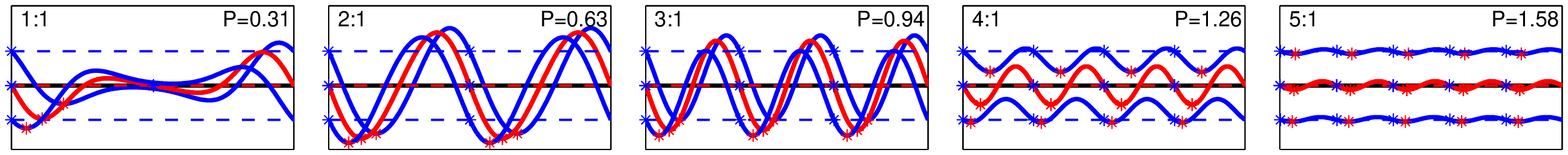}} \vspace{-.25cm}
\caption{\footnotesize Interaction functions $\Lambda_{v}^{\uNM}$ for fast subharmonic entrainment controls $v_f$ for the Hodgkin-Huxley neuron model with frequency detuning of {\color{red}$0\%$} and {\color{blue}$\pm5\%$}, where the domain and range in each plot is $[0,2\pi]$ and $[-0.04,0.05]$, respectively.  The entrainment ratio and control energy are indicated, and the black line indicates the $x$-axis.  As $N$ grows large, the interaction function converges to small variation about $\Delta\w$, which is indicated in each case by a dashed line.  Red and blue stars mark the phases $\{\vphi^-_i\}$ and stable fixed points $\{\vphi_i^*\}$ of \eqref{sys5}, each of which occur $N$ times. } \vspace{.25cm} \label{s5f2_fif}
\end{figure}

Substituting \eqref{sol:cand} into the equation \eqref{dlnmdf}, and recalling that $\bt{Y_\vphi^{\,\uNM} Y^{\,\uNM}}=0$, yields
\begin{align} \label{fenm4}
\ddx{\vphi}\Lambda_{v_f}^{\uNM}(\vphi^*) & = \frac{1}{2\lambda}\bq{\bt{Y_\vphi^{\,\uNM}(\eta,\vphi^*)Y_\vphi^{\,\uNM}(\eta,\vphi^*)} - \mu \bt{Y_\vphi^{\,\uNM}(\eta,\vphi^*)Y^{\,\uNM}(\eta,\vphi^*)}  } \notag
%\frac{1}{2\pi M}\sum_{\ell=0}^{M-1}\int_0^{2\pi} Z'\bp{\xi + \vphi^*} \cdot \nonumber \\ & \qquad \frac{1}{2\lambda}\bq{Y_\vphi^{\,\uNM}\bp{\frac{N}{M}[2\pi\ell+\xi],\vphi^*} - \mu Y^{\,\uNM}\bp{\frac{N}{M}[2\pi\ell+\xi],\vphi^*} }\rd \xi \notag \\
%&=  \frac{1}{4\pi N \lambda}\sum_{j=0}^{N-1}\int_0^{2\pi} Z'(\xi + \vphi^*)Z'\bp{\xi + \frac{M}{N}2\pi j+\vphi^*}\rd \xi  \notag
\\ &= \frac{1}{2\lambda}\bt{Y_\vphi^{\,\uNM}Y_\vphi^{\,\uNM}} = \frac{1}{2\lambda}S_0^{\,\uNM}.
\end{align}
In particular, $S_0^{\,\uNM}>0$, so we choose $\lambda<0$ when solving (\ref{lfasteq}) for $\lambda$ in order for the expression \eqref{fenm4} to be negative in order for the objective in (\ref{prob:1}) to be maximized.
It follows that the optimal waveform and multiplier can be obtained from \eqref{sol:cand}, \eqref{mufasteq} and \eqref{lfasteq} as
\begin{align} \label{sol:opt}
v_f(\eta)=\frac{Y_\vphi^{\,\uNM}(\eta,\vphi^*)}{2\lambda}- \frac{\Delta \w Y^{\,\uNM}(\eta,\vphi^*)}{V_0^{\,\uNM}} , \quad \lambda = -\frac{1}{2}\sqrt{\frac{S_0^{\,\uNM}}{P-\frac{(\Delta\w)^2}{V_0^{\,\uNM}}}},
\end{align}
%we disregard the phase shift $\vphi^*$, because entrainment is asymptotic.
where phase-locking occurs fastest when the oscillator is in the neighborhood of the phase $\psi(t)=\vphi^*$ at the start of entrainment with $v_f$.   For zero frequency detuning, the optimal waveform is a re-scaling of $Y_\vphi^{\,\uNM}$, which is a sum of shifted derivatives of the PRC function.  As $|\Delta\w|$ increases, $v$ continuously transforms towards a rescaling of $Y^{\,\uNM}$, which is the minimum energy waveform for subharmonic entrainment, as derived in the previous section.  This transition reflects the conceptual trade-off between the fast entrainment objective \eqref{prob:1} and frequency control constraint \eqref{eq:con2}, which can be satisfied only when $P>(\Delta\w)^2/V_0^{\,\uNM}$, as shown in \eqref{minpvpow}.  When $N=M=1$, these results reduce to the harmonic (1:1) case \cite{zlotnik13prl}. Figure \ref{s5f1_fk} shows the fast subharmonic entrainment controls for the Hodgkin-Huxley model for values of $N,\,M=1,\ldots,5$ at detunings between $-5\%$ and $5\%$, and Figure \ref{s5f2_fif} shows several corresponding interaction functions.  It is important to note that the choice of control energy $P$ significantly impacts the control waveform in the case of non-zero detuning, as seen in the panels on the diagonal of Figure \ref{s5f1_fk}.

\section{Minimum energy subharmonic entrainment of oscillator ensembles} \label{secnmens}

In practice, biological systems exhibit variation in parameters that characterize the system dynamics, which must be taken into account when designing optimal entrainment controls.  We approach this issue by modeling an ensemble of systems as a collection of phase models with a common PRC that is derived using a nominal parameter set, and the frequencies span the range of natural frequencies resulting from phase model reduction of systems with parameters in a specified range.  A justification of this approach and a sensitivity analysis is provided in Section 5 in \cite{zlotnik12jne}.  The following extension of this modeling and control technique to subharmonic ($N$:$M$) entrainment parallels our previous work \cite{zlotnik12jne} by incorporating the theory derived in Section \ref{secent}, and contains a more rigorous optimality proof and additional generalizations.

Specifically, we consider a collection of systems $\dot{x}=f(x,u,p)$ where $x\in\bR^n$ is the state, $u\in\bR$ is a scalar control, and $p\in\cD\subset\bR^d$ is a vector of constant parameters varying on a hypercube $\cD$ containing a nominal parameter vector $q$.  Each system can be reduced to a scalar phase model $\dot{\theta}=\w(p)+Z(\theta,p)u$, where the natural frequency and PRC depend on the parameter vector $p$.   In order to design a control that entrains the ensemble for all $p\in\cD$, we approximate it by $\{\dot{\theta}=\w(p)+Z(\theta,q)u\,\,:\,\,p\in\cD\}$, where $Z(\theta,q)$ is the nominal PRC.

\begin{figure}[t]
\hspace{29pt} \centerline {\includegraphics[width=.75\linewidth]{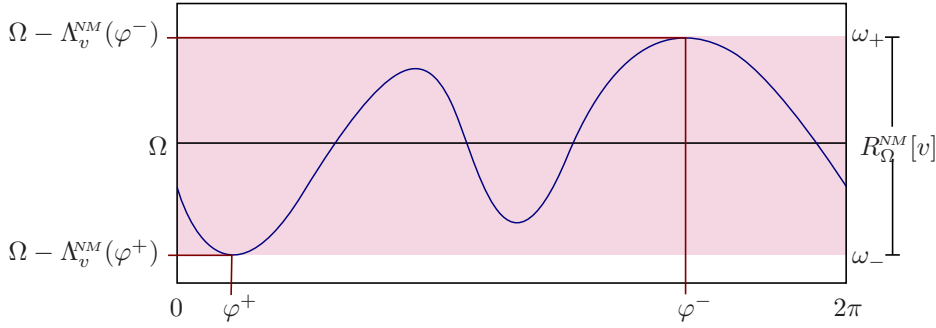}}
\begin{picture}(0,0)
\put(18,109){$\Omega-\Lambda_v^{\uNM}(\vphi^-)$} \put(18,26){$\Omega-\Lambda_v^{\uNM}(\vphi^+)$}
\put(337,109){$\w_+$} \put(337,25){$\w_-$}
\put(99,5){$\vphi^+$} \put(271,5){$\vphi^-$}
\put(79,5){$0$} \put(330,5){$2\pi$} \put(71,66){$\Omega$} \put(340,66){$R_\O^{\,\uNM}[v]$}
\end{picture}
\caption{\footnotesize This graphic illustrates the constraints \eqref{aconst2}.  The curve shown is $\O-\Lambda_v^{\uNM}(\vphi)$, and the frequency locking region $R_\O^{\,\uNM}[v]=[\w_-,\w_+]$ is indicated by pink shading. When $[\w_1,\w_2]\subseteq R_\O^{\,\uNM}[v]$, then the collection \eqref{en_sys_1} is entrained to $\Omega$. } \vspace{.25cm} \label{s6f1_inter1}
\end{figure}

Our strategy is to derive a minimum energy periodic control signal $u(t)=v(\O_f t)$ that guarantees entrainment for each system in the ensemble of oscillators
\begin{align}
\cF=\{\dot{\psi}=\w+Z(\psi)u\,\,:\,\,\w\in[\w_1,\w_2]\}
\label{en_sys_1}
\end{align}
to a frequency $\O$, where the target and forcing frequencies satisfy $M\Omega_f=N\Omega$.  We approach the subharmonic entrainment of oscillator ensembles by applying the theory in Section \ref{secent} to the derivation of optimal ensemble controls in Section 4 of \cite{zlotnik12jne}.    We call the range of frequencies that are entrained by the control $v$ applied at the frequency $\O_f$ with subharmonic ratio $N$:$M$ the subharmonic locking range $R_\O^{\,\uNM}[v]=[\w_-,\w_+]$, and when $[\w_1,\w_2]\subseteq R_\O^{\,\uNM}[v]$ we say that the ensemble $\cF$ is entrained.   This requirement results in two constraints, which can be visualized with the help of Figure \ref{s6f1_inter1}, of the form
\begin{align} \label{aconst2}
\begin{array}{rcrcrcrcr}
\Delta \w_+ &\teQ & \w_+-\O &=&-\Lambda_v^{\uNM}(\vphi^-) &\geq& \w_2-\O &\teQ & \Delta \w_2,  \\
\Delta \w_- &\teQ &\w_--\O &=&-\Lambda_v^{\uNM}(\vphi^+) & \leq &  \w_1-\O & \teQ & \Delta \w_1.
\end{array}
\end{align}

\begin{figure}[t]
\vspace{.5cm}
\hspace{10pt} \centerline {\includegraphics[width=.92\linewidth]{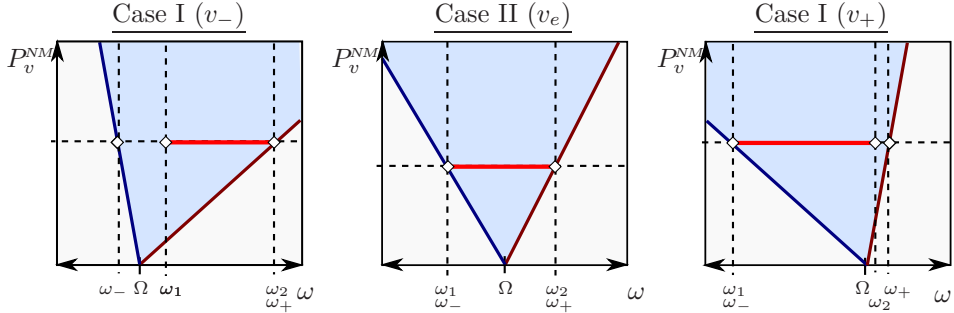}}
\begin{picture}(0,0)
\put(10,91){$P_v^{\,\uNM}$} \put(45,5){$\scriptstyle \w_{\scriptscriptstyle -}$}  \put(58,5){$\scriptstyle\O$} \put(68,5){$\scriptstyle\w_1$} \put(68,5){$\scriptstyle\w_1$} \put(108,5){$\scriptstyle \w_2$} \put(108,0){$\scriptstyle \w_+$} \put(120,3){$\w$} \put(50,108){\underline{Case I ($v_-$)}}
\put(133,91){$P_v^{\,\uNM}$} \put(172,5){$\scriptstyle\w_1$} \put(172,0){$\scriptstyle\w_-$} \put(196,5){$\scriptstyle\O$} \put(214,5){$\scriptstyle\w_2$} \put(214,0){$\scriptstyle\w_+$}  \put(245,3){$\w$} \put(172,108){\underline{Case II ($v_e$)}}
\put(256,91){$P_v^{\,\uNM}$} \put(281,5){$\scriptstyle\w_1$} \put(281,0){$\scriptstyle\w_-$} \put(331,5){$\scriptstyle\O$} \put(336,0){$\scriptstyle\w_2$} \put(342,5){$\scriptstyle\w_+$} \put(361,2){$\w$} \put(293,108){\underline{Case I ($v_+$)}}
\end{picture}
\caption{\footnotesize Illustration of ensemble Arnold tongues for Case I (both $v_-$ and $v_+$) and Case II controls, which are used when constraints on $\Lambda_v^{NM}(\vphi^-)$ and $\Lambda_v^{NM}(\vphi^+)$ are active.  Note that one-sided tongues as in Figure \ref{s3f2_arnold} can occur in either case, depending on $N$, $M$, and the PRC $Z$, as shown in Section \ref{sechh}.  The range $[\w_1,\w_2]$ of frequencies in the ensemble $\cF$ is marked with a {\color{red}red} bar.}  \vspace{.25cm} \label{s6f2_arnold}
\end{figure}

The objective of minimizing control energy $\bt{v^2}$ given the constraints (\ref{aconst2}) gives rise to the optimization problem
\begin{equation} \label{aop2}
\begin{array}{rl}
\min & \cJ[v] = \bt{v^2}, \quad v\in\cP\medskip\\
\st &\,\,\,\, \Delta \w_2+\Lambda_v^{\uNM}(\vphi^-) \leq 0, \medskip \\
 & -\Delta \w_1-\Lambda_v^{\uNM}(\vphi^+) \leq 0.
\end{array}
\end{equation}
We refer to the event that one of $v_+$ (when $\w_2<\O$) or $v_-$ (when $\w_1>\O)$ in \eqref{solnm1} can solve the problem \eqref{aop2} as Case I.   Understanding the Arnold tongues that characterize subharmonic entrainment of ensembles in the form of $\cF$, as illustrated in Figure \ref{s6f2_arnold}, will clarify the conditions when \eqref{solnm1} is optimal, and when another class of solutions, which we call Case II, is superior.  We derive this condition, which depends on the ensemble parameters $\w_1$ and $\w_2$ and $Z$ as well as the target frequency $\O$ and the subharmonic ratio $N$:$M$.

In contrast to the description in Section \ref{secent} of the Arnold tongue associated with a single oscillator and a given waveform, given an ensemble $\cF$ we are interested in the relationship between the locking range $R_\O^{\,\uNM}[v]$ and the RMS control power.  Therefore, we define the ensemble Arnold tongue as the set of pairs $(\omega,P_v^{NM})\in\bR^2$ that result in entrainment of an oscillator in $\cF$ with natural frequency $\w$ to a frequency $\O_f$ at a subharmonic $N$:$M$ ratio using the waveform $v$, where $w$ is the natural frequency of the oscillator.  The equation \eqref{arnold_tg1} is modified to
$\w-\O+P_v^{\,\uNM}(\w)\cdot\Lambda_{\tw{v}}^{\uNM}(\vphi)=0$, where the left and right boundaries of the Arnold tongue are approximated by solving for $P_v^{\,\uNM}(\w)$ as a function of $\w$ and substituting $\vphi=\vphi^-$ and $\vphi=\vphi^+$, respectively.  This yields
\begin{align} \label{ensatng1}
P_v^{\,\uNM}(\w)=\left\{\begin{array}{ll} (\O-\w)/\Lambda_{\tw{v}}^{\uNM}(\vphi^-), & \rM{right} \\  (\O-\w)/\Lambda_{\tw{v}}^{\uNM}(\vphi^+), & \rM{left}, \end{array}\right.
\end{align}
as a linear estimate of the ensemble Arnold tongue boundary, where $\tw{v}=v/\sqrt{\bt{v^2}}$ is the unit power normalization of $v$ as before.  Illustrations of Arnold tongues for the two possible cases are illustrated in Figure \ref{s6f2_arnold}.  The notion of Arnold tongues guides our derivation in the following subsections of the possible optimal control solutions and criteria for optimality of these different cases.

\subsection{Case I: Solution \eqref{solnm1} is optimal for subharmonic ensemble entrainment} \label{enscase1}

To derive the conditions when \eqref{solnm1} is optimal, we focus on the entrainment of $\cF$ to  a frequency $\O\in[\w_1,\w_2]$ using $v_-$ when $\Delta\w_+=\Delta\w_2>-\Delta\w_1$ (i.e., when $\w_2$ is further from $\O$ than $\w_1$) while noting that the case where $v_+$ is optimal for $\Delta\w_2<-\Delta\w_1=-\Delta\w_-$ (i.e., when $\w_1$ is further from $\O$ than $\w_2$) is symmetric.  Because $\w_2$ is the  natural frequency in the ensemble farthest from $\O$, we use $\Delta\w=\Delta\w_2$ in \eqref{solnm1}, and then check whether $[\w_1,\w_2]\subseteq R_\Omega^{\,\uNM}$.  Then the first constraint in (\ref{aop2}) is active, which yields
\begin{align} \label{avminent1}
-\Delta\w_+=\Lambda_v^{\uNM}(\vphi^-)=-\Delta\w_2,
\end{align}
so that $\w_+=\w_2$ is the upper bound on the locking range $R_\O^{\,\uNM}[v]$, as desired. It remains to determine $\Lambda_v^{\uNM}(\vphi^+)=\O-\w_-$, from which we obtain the lower bound $\w_-$ on $R_\O^{\,\uNM}[v]$.  Using the expression \eqref{lnm2b} for $\Lambda_v^{\uNM}(\vphi)$ together with the solution for $v_-$ in \eqref{solnm1} using $\Delta\w=\Delta\w_2$, we find that
\begin{align}
\Lambda_{v_-}^{\uNM}(\vphi)  =  \bt{Y^{\,\uNM}(\eta,\vphi)v_-(\eta)} =
-\frac{\Delta\w_2}{V_0^{\,\uNM}}V^{\,\uNM}(\vphi-\vphi^-).
\end{align}
Observe that $\Lambda_v^{\uNM}(\vphi)$ is maximized when $V^{\,\uNM}(\vphi-\vphi^-)$ is minimized, and hence to find $\Lambda_v^{\uNM}(\vphi^+)$ it suffices to find the minimum value $V_*^{\,\uNM}$ of $V^{\,\uNM}$.
It follows that
\begin{align} \label{abound1}
\Lambda_{v_-}^{\uNM}(\vphi^+)=-\frac{\Delta\w_2}{V_0^{\,\uNM}}V_*^{\,\uNM},
\end{align}
and the lower bound of $R_\O^{\,\uNM}[v]$ is $\w_-=\O-\Lambda_v^{\uNM}(\vphi^+)$. If $\w_-<\w_1$, then $[\w_1,\w_2]\subseteq R_\O^{\,\uNM}[v]$, hence the control $v_-$ in (\ref{solnm1}), with $\Delta\w = \w_2-\O$, is the minimum energy solution to problem (\ref{aop2}), and entrains $\cF$ to the frequency $\O$.

Therefore to determine whether the problem is optimally solved by $v_-$, the decision criterion is obtained by combining the definition $\Lambda_v^{\uNM}(\vphi^+)=-\Delta\w_-$  with \eqref{avminent1} and \eqref{abound1} to yield the boundary estimate $\Delta\w_-=\Delta\w_+ V_*^{\,\uNM}/V_0^{\,\uNM}$. Thus if the relation
\begin{align} \label{aomegrat1}
\Delta\w_1>\Delta\w_-=\frac{\Delta\w_2}{V_0^{\,\uNM}}V_*^{\,\uNM}
\end{align}
is satisfied, then $v_-$ with $\Delta\w=\Delta\w_2$ will be optimal.  The derivation of the condition when $v_+$ is optimal is symmetric, and results in a boundary estimate $\Delta\w_+=\Delta\w_1V_*^{\,\uNM}/V_0^{\,\uNM}$.  It follows that if the condition
\begin{align} \label{aomegrat2}
\Delta\w_2<\Delta\w_+=\frac{\Delta\w_1}{V_0^{\,\uNM}}V_*^{\,\uNM}
\end{align}
holds, then the control $v_+$ with $\Delta\w=\Delta\w_1$ is optimal for entraining $\cF$ to the frequency $\O$. When neither \eqref{aomegrat1} or \eqref{aomegrat2} holds, then neither $v_-$ or $v_+$ in \eqref{solnm1} is the solution to \eqref{aop2}.  In the following subsection, we derive the optimal solution for that case.

\subsection{Case II: Neither of conditions \eqref{aomegrat1} and \eqref{aomegrat2} are satisfied} \label{enscase2}

In Case I above, when \eqref{solnm1} is optimal for entraining the ensemble \eqref{en_sys_1}, only one of the constraints in problem \eqref{aop2} is active.  When neither of the conditions \eqref{aomegrat1} and \eqref{aomegrat2} is satisfied, the solution to problem \eqref{aop2} occurs when both constraints are active.  To derive this solution, we adjoin the constraints in \eqref{aop2} to the minimum energy objective function using multipliers $\mu_-$ and $\mu_+$, which gives rise to the cost functional
\begin{align} \label{aop3}
\cJ[v] & =  \bt{v^2}-\mu_-(\Delta \w_2+\Lambda_v^{\uNM}(\vphi^-)) -\mu_+(-\Delta \w_1-\Lambda_v^{\uNM}(\vphi^+)) \notag\\
& =  \frac{1}{2\pi}\int_0^{2\pi} \Bigl(v(\eta)[v(\eta) - \mu_-Y^{\,\uNM}(\eta,\vphi^-) +\mu_+Y^{\,\uNM}(\eta,\vphi^+) ] \notag \\
&  \qquad -\mu_-\Delta\w_2+\mu_+\Delta\w_1 \Bigr)\rd \eta,
\end{align}
where we have used the expression \eqref{lnm2b} for $\Lambda_v^{\uNM}$.
Solving the Euler-Lagrange equation yields
\begin{align} \label{asol3}
v_e(\theta)=\dS-\frac{1}{2}[ \mu_+Y^{\,\uNM}(\eta,\vphi^+)-\mu_-Y^{\,\uNM}(\eta,\vphi^-)],
\end{align}
which we substitute back into problem (\ref{aop2}) to obtain
\begin{eqnarray} \label{aeq3}
\bt{v_e^2}&=&\frac{1}{4}\langle\bp{\mu_+Y^{\,\uNM}(\eta,\vphi^+)-\mu_-Y^{\,\uNM}(\eta,\vphi^-)}^2\rangle  %\\ & = & \frac{1}{4}\mu_+^2 \bt{Z^2} -\half\mu_+\mu_-\bt{Z(\theta+\vphi^+)Z(\theta+\vphi^-)} + \frac{1}{4}\mu_-^2 \bt{Z^2} \notag
\\ & = & \frac{1}{4}(\mu_+^2+\mu_-^2)V_0^{\,\uNM} - \half \mu_+\mu_-V^{\,\uNM}(\Delta\vphi), \notag\\\notag\\
\Lambda_{v_e}^{\uNM}(\vphi^+) & = &  -\half\mu_+V_0^{\,\uNM} + \half\mu_-V^{\,\uNM}(\Delta\vphi), \label{aeq4} \\
\Lambda_{v_e}^{\uNM}(\vphi^-) & = &  \half \mu_-V_0^{\,\uNM}-\half\mu_+V^{\,\uNM}(\Delta\vphi), \label{aeq5}
\end{eqnarray}\\
where $\Delta\vphi=\vphi^+-\vphi^-$ is the range spanned by $\Lambda_v^{\uNM}(\vphi)$ for $\vphi\in[0,2\pi]$. Using the expressions (\ref{aeq3}), (\ref{aeq4}), and (\ref{aeq5}) transforms the functional optimization problem (\ref{aop2}) into a nonlinear programming problem in the variables $\mu_-$, $\mu_+$, and $V^{\,\uNM}(\Delta\vphi)$, given by
\begin{equation} \label{aop4}
\begin{array}{rr}
\min & \cJ[\mu_-,\mu_+,V(\Delta\vphi)] = \dS\tfrac{1}{4}(\mu_+^2+\mu_-^2)V_0^{\,\uNM} - \tfrac{1}{2} \mu_+\mu_-V^{\,\uNM}(\Delta\vphi) \medskip \\
\st &\dS\Delta \w_2+\tfrac{1}{2} \mu_-V_0^{\,\uNM}-\tfrac{1}{2}\mu_+V^{\,\uNM}(\Delta\vphi)\leq 0, \medskip \\
 & \dS-\Delta \w_1+\tfrac{1}{2}\mu_+V_0^{\,\uNM} - \tfrac{1}{2}\mu_-V^{\,\uNM}(\Delta\vphi)\leq 0.
\end{array}
\end{equation}
We focus in Case II on optimal solutions to problem \eqref{aop4} for which both constraints are active.  Indeed, when Case I is in effect, one of conditions \eqref{aomegrat1} or \eqref{aomegrat2} is satisfied, so that  $\mu_+=0$ or $\mu_-=0$, and problem \eqref{aop4} is reduced to problem \eqref{op1} with $\lambda=\mu_-$ or $\lambda=-\mu_+$, respectively.  Otherwise, both constraints in problem \eqref{aop4} are active, with multipliers given by
\begin{equation} \label{amult1}
\begin{array}{rcl}
\mu_+ & =&\dS\frac{2(\Delta\w_1V_0^{\,\uNM}-\Delta\w_2 V^{\,\uNM}(\Delta\vphi))}{(V_0^{\,\uNM}-V^{\,\uNM}(\Delta\vphi))(V_0^{\,\uNM}+V^{\,\uNM}(\Delta\vphi))}, \medskip \\  \mu_- & =&\dS\frac{2(\Delta\w_1V^{\,\uNM}(\Delta\vphi)-\Delta\w_2V_0^{\,\uNM}) }{(V_0^{\,\uNM}-V^{\,\uNM}(\Delta\vphi))(V_0^{\,\uNM}+V^{\,\uNM}(\Delta\vphi))}.
\end{array}
\end{equation}
For these multipliers, the objective in problem (\ref{aop4}) is reduced to function of $\beta=V^{\,\uNM}(\Delta\vphi)$ given by
\begin{eqnarray} \label{acost1}
\cJ[\beta] &=& \dS\frac{(\Delta\w_1V_0^{\,\uNM}-\Delta\w_2 \beta )^2+ (\Delta\w_1 \beta -\Delta\w_2V_0^{\,\uNM})^2}{(V_0^{\,\uNM}- \beta )^2(V_0^{\,\uNM}+ \beta )^2} V_0^{\,\uNM}  \notag\\ & & - \dS\frac{2(\Delta\w_1 V_0^{\,\uNM}-\Delta\w_2 \beta )(\Delta\w_1 \beta -\Delta\w_2 V_0^{\,\uNM}) \beta }{(V_0^{\,\uNM}- \beta )^2( V_0^{\,\uNM}+ \beta )^2}.
\end{eqnarray}
Differentiating the cost (\ref{acost1}) with respect to $\beta$ results in\\
\begin{align} \label{adiff1}
\dxdy{\cJ[\beta ]}{\beta } = -2\frac{(V_0^{\,\uNM}\Delta\w_1-\beta\Delta\w_2)(V_0^{\,\uNM}\Delta\w_2-\beta\Delta\w_1)}{(V_0^{\,\uNM}-\beta )^2(V_0^{\,\uNM}+\beta )^2}.
\end{align}
Recall that because neither of the conditions \eqref{aomegrat1} or \eqref{aomegrat2} holds, then
\begin{align} \label{notaomegrat}
V_0^{\,\uNM}\Delta\w_1-V_*^{\,\uNM}\Delta\w_2<0 \quad \anD \quad V_0^{\,\uNM}\Delta\w_2-V_*^{\,\uNM}\Delta\w_1>0.
\end{align}
We are restricted to $\beta=V^{\,\uNM}(\Delta\vphi)\in[V_*^{\,\uNM},V_0^{\,\uNM}]$, so we write $V_*^{\,\uNM}-\beta = p(V_*^{\,\uNM}-V_0^{\,\uNM})$ for some $p\in[0,1]$.  In addition, the fact that $\w_1<\w_2$ results in $\Delta\w_1-\Delta\w_2<0$ and $\Delta\w_2-\Delta\w_1>0$. Therefore the quantities in the numerator of \eqref{adiff1} satisfy
\begin{align} \label{adiffnum1}
V_0^{\,\uNM}\Delta\w_1-\beta\Delta\w_2 & = V_0^{\,\uNM}\Delta\w_1-V_*^{\,\uNM}\Delta\w_2 + (V_*^{\,\uNM}-\beta)\Delta\w_2  \\  & = V_0^{\,\uNM}\Delta\w_1-V_*^{\,\uNM}\Delta\w_2 + p(V_*^{\,\uNM} - V_0^{\,\uNM})\Delta\w_2 \notag \\ & = p(\Delta\w_1-\Delta\w_2)V_0^{\,\uNM} + (1-p)(V_0^{\,\uNM}\Delta\w_1-V_*^{\,\uNM}\Delta\w_2) < 0, \notag
\end{align}
\begin{align} \label{adiffnum2}
V_0^{\,\uNM}\Delta\w_2-\beta\Delta\w_1  & = V_0^{\,\uNM}\Delta\w_2-V_*^{\,\uNM}\Delta\w_1 + (V_*^{\,\uNM} - \beta)\Delta\w_1  \\ & =  V_0^{\,\uNM}\Delta\w_2-V_*^{\,\uNM}\Delta\w_1 + p(V_*^{\,\uNM} - V_0^{\,\uNM})\Delta\w_1 \notag \\ & = p(\Delta\w_2-\Delta\w_1)V_0^{\,\uNM}+ (1-p)(V_0^{\,\uNM}\Delta\w_2-V_*^{\,\uNM}\Delta\w_1)>0. \notag
\end{align}
The relations \eqref{adiffnum1} and \eqref{adiffnum2} imply that \eqref{adiff1} is positive for all values of $\beta=V^{\,\uNM}(\Delta\vphi)\in[V_*^{\,\uNM},V_0^{\,\uNM}]$, so that the cost (\ref{acost1}) increases when $V^{\,\uNM}(\Delta\vphi)$ does.  Therefore the objective (\ref{acost1}) is minimized when $V^{\,\uNM}(\Delta\vphi)$ is, which occurs when $V^{\,\uNM}(\Delta\vphi)=V_*^{\,\uNM}$.  Therefore the problem (\ref{aop4}) is solved when \begin{align} \label{magang1}
\Delta\vphi=\vphi_*^{\uNM}=\underset{\vphi\in[0,2\pi]}{\rM{argmin}}V^{\,\uNM}(\vphi),
\end{align}
and the multipliers are as in (\ref{amult1}).  The locking range for this control is then exactly $R_\O^{\,\uNM}[v_e]=[\w_1,\w_2]$, which satisfies the entrainment constraints (\ref{aconst2}).

By combining the results in Sections \ref{enscase1} and \ref{enscase2} for Cases I and II we can completely characterize the minimum energy control that entrains the ensemble $\cF$ in \eqref{en_sys_1} to a target frequency $\Omega$ with subharmonic ratio $N$:$M$.  This full solution is
{\setstretch{1.2}
\begin{alignat}{1} \label{asoltot}
v_e(\eta)& = \left\{\begin{array}{ll} \dS -\frac{\Delta \w_1}{V_0^{\,\uNM}} Y^{\,\uNM}(\eta,\vphi^+) &  \, \text{if}\quad \dS \Delta\w_2<\frac{\Delta\w_1}{V_0^{\,\uNM}}V_*^{\,\uNM}, \\  \\
\dS {\def\arraystretch{1.7}\begin{array}{l} \dS
\frac{(\Delta\w_2 V_*^{\,\uNM}-\Delta\w_1V_0^{\,\uNM})}{(V_0^{\,\uNM}-V_*^{\,\uNM})(V_0^{\,\uNM}+V_*^{\,\uNM})}Y^{\,\uNM}(\eta,\vphi_*^{\uNM}) \\  \dS \qquad + \,\, \frac{(\Delta\w_1V_*^{\,\uNM}-\Delta\w_2V_0^{\,\uNM})}{(V_0^{\,\uNM}-V_*^{\,\uNM})(V_0^{\,\uNM}+V_*^{\,\uNM})}Y^{\,\uNM}(\eta,0)  \end{array} }
&\, \text{if}\quad \left\{\begin{array}{l} \dS \Delta\w_1<\frac{\Delta\w_2}{V_0^{\,\uNM}}V_*^{\,\uNM} \\ {\text{\footnotesize and}} \\ \dS \Delta\w_2>\frac{\Delta\w_1}{V_0^{\,\uNM}}V_*^{\,\uNM} \end{array}\right. \\ \\\dS
-\frac{\Delta \w_2}{V_0^{\,\uNM}}Y^{\,\uNM}(\eta,\vphi^-) &  \, \text{if}\quad  \dS \Delta\w_1>\frac{\Delta\w_2}{V_0^{\,\uNM}}V_*^{\,\uNM}. \end{array}\right.
\end{alignat}}
%Note that it is sufficient in practice to maintain a phase shift of $\Delta\vphi=\vphi^+-\vphi^-$ between the weighted $Y^{\,\uNM}$ functions, because entrainment is an asymptotic process and the initial phase ambiguity can be omitted.
Finally, the energy of $v_e$, which is the minimum value of the objective \eqref{aop3}, simplifies to
\begin{alignat}{1} \label{asoltote}
\bt{v_e^2}& = \left\{\begin{array}{ll} \dS \frac{(\Delta\w_1)^2}{V_0^{\,\uNM}} &  \quad \text{if}\quad \dS \Delta\w_2<\frac{\Delta\w_1}{V_0^{\,\uNM}}V_*^{\,\uNM}, \\
\dS \frac{(\Delta\w_1^2+\Delta\w_2^2)V_0^{\,\uNM}-2\Delta\w_1\Delta\w_2 V_*^{\,\uNM}}{(V_0^{\,\uNM}-V_*^{\,\uNM})(V_0^{\,\uNM}+V_*^{\,\uNM})}
&\quad \text{if}\quad \left\{\begin{array}{l} \dS \Delta\w_1<\frac{\Delta\w_2}{V_0^{\,\uNM}}V_*^{\,\uNM} \\ {\text{\footnotesize and}} \\ \dS \Delta\w_2>\frac{\Delta\w_1}{V_0^{\,\uNM}}V_*^{\,\uNM} \end{array}\right. \\ \dS
\frac{(\Delta\w_2)^2}{V_0^{\,\uNM}}  &  \quad \text{if}\quad  \dS \Delta\w_1>\frac{\Delta\w_2}{V_0^{\,\uNM}}V_*^{\,\uNM}. \end{array}\right.
\end{alignat}

We have shown that the minimum energy periodic control $u(t)=v(\frac{N}{M}\O t)$ that achieves subharmonic entrainment of an ensemble of oscillators \eqref{en_sys_1} to a target frequency $\O$ is an appropriately weighted sum of shifted functions $Y^{\,\uNM}$, as given in \eqref{ynmdf}, where $\eta=\frac{N}{M}\O t=\O_f t$ is the forcing phase. When $N=M=1$, these results reduce to the optimal solution for the harmonic (1:1) case \cite{zlotnik12jne}.  Figure \ref{s6f3_wk} shows the minimum energy subharmonic controls for ensembles of Hodgkin-Huxley neurons for $N,M=1,\ldots,5$ and several ranges of $[\w_1,\w_2]$.  We have also presented generalized criteria for using the two derived classes of optimal controls, which can be applied to systems with Type I (strictly positive) and Type II PRCs, while the derivation in our previous work on $1$:$1$ entrainment required $Q_*<0$.  It is important to note that using $\Omega=\half(\w_1+\w_2)$ allows the ensemble to be entrained with a minimum control energy, as in the case of $[\w_1,\w_2]=${\color{rb3} $[0.95\w,1.05\w]$}.  One can then consider a dual objective of maximizing the locking range of entrainment given a fixed power, as described in the following section.

\section{Maximum locking range for subharmonic entrainment} \label{secmaxr}

In some applications, the frequency range $[\w_1,\w_2]$ of the oscillators in $\cF$ in \eqref{en_sys_1}  is not known, or it is desirable to entrain the largest collection of oscillators with similar PRC but uncertain frequency.  For such cases, we seek a control $v$ that maximizes the locking range of entrainment $R_\O^{\,\uNM}[v]$ for a fixed control energy $\bt{v^2}=P$.  Because $R_\O^{\,\uNM}[v]=[\w_-,\w_+]$, we wish to maximize $\w_+-\w_-=\Delta\w_+-\Delta\w_-=\Lambda_{v}^{\uNM}(\vphi^+)-\Lambda_{v}^{\uNM}(\vphi^-)$, where the latter equality is due to the constraints \ref{aconst2}.  The resulting optimization problem can be formulated as
\begin{align} \label{prob:wide}
	\max_{v\in\cP} \quad & \cJ[v] = \Lambda_{v}^{\uNM}(\vphi^+)-\Lambda_{v}^{\uNM}(\vphi^-) \\
	\label{eq:conw1}
	\rM{s.t.} \quad & \bt{v^2} = P.
\end{align}
By adjoining the constraint to \eqref{eq:conw1} to the objective using a multiplier $\lambda$, we obtain a cost functional given by
\begin{align} \label{prob:wide2}
	\cJ[v] &= \Lambda_{v}^{\uNM}(\vphi^+)-\Lambda_{v}^{\uNM}(\vphi^-) -\lambda(\bt{v^2} -P) \notag \\ & = \bt{Y^{\,\uNM}(\eta,\vphi^+)v(\eta)} - \bt{Y^{\,\uNM}(\eta,\vphi^-)v(\eta)} -\lambda(\bt{v^2} -P) \notag \\ & = \frac{1}{2\pi}\int_0^{2\pi} \bp{ v(\eta)\bq{Y^{\,\uNM}(\eta,\vphi^+) - Y^{\,\uNM}(\eta,\vphi^-) -\lambda v(\eta)} +\lambda P} \rd\eta.
\end{align}
Solving the Euler-Lagrange equation yields a candidate solution in the form
\begin{align} \label{sol:candw}
 v_r(\eta)=\frac{1}{2\lambda}\bq{Y^{\,\uNM}(\eta,\vphi^+)-Y^{\,\uNM}(\eta,\vphi^-)}.
\end{align}
By applying \eqref{lnmv1}, the interaction function is shown to be
\begin{align} \label{sol:candw2}
\Lambda_{v_r}^{\uNM}(\vphi)=\frac{1}{2\lambda}[V^{\,\uNM}(\vphi-\vphi^+)-V^{\,\uNM}(\vphi-\vphi^-)],
\end{align}
so the objective \eqref{prob:wide} is given by
\begin{align} \label{objw2}
\Lambda_{v_r}^{\uNM}(\vphi^+)-\Lambda_{v_r}^{\uNM}(\vphi^-) = \frac{1}{\lambda}[V^{\,\uNM}(0)-V^{\,\uNM}(\Delta\vphi)].
\end{align}
\begin{figure}[t]
\centerline {\includegraphics[width=\linewidth]{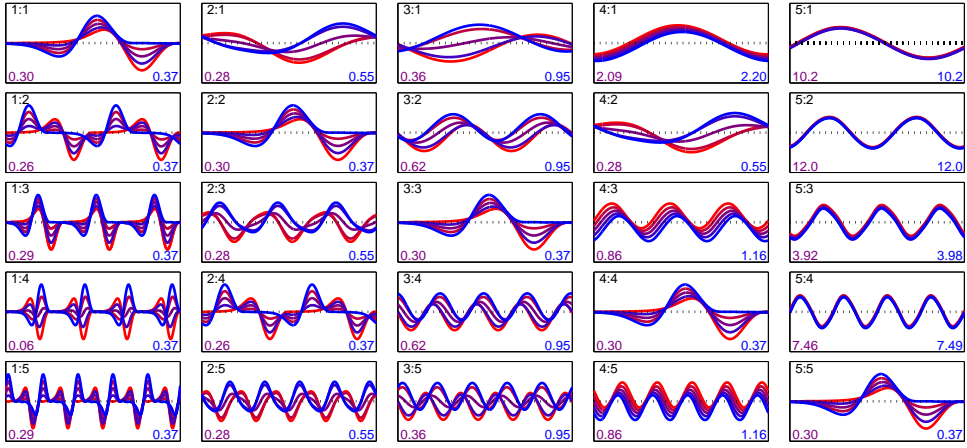}} \vspace{-.25cm}
\caption{\footnotesize Minimum energy subharmonic entrainment controls for Hodgkin-Huxley neuron ensembles with frequency ranges {\color{rb1} $[0.925\w,1.025\w]$}, {\color{rb2} $[0.9375\w,1.0375\w]$}, {\color{rb3} $[0.95\w,1.05\w]$}, {\color{rb4} $[0.9625\w,1.0625\w]$}, and {\color{rb5} $[0.975\w,1.075\w]$}, and target frequency $\O=\w$ equal to the nominal natural frequency in Appendix \ref{aphh}.  In each panel, the controls are rescaled so that the greatest energy waveform has unit energy, and the domain and range in each plot is $[0,2\pi]$ and $[-3.7,3.7]$, respectively.  The entrainment ratio is indicated at the top, while the lowest energy (for {\color{rb3}$v_r$}) and highest energy (for {\color{rb1}$v_-$} or {\color{rb5}$v_+$}) are shown at bottom left and right, respectively.} \vspace{.25cm} \label{s6f3_wk}
\end{figure}
\begin{figure}[t]
\centerline {\includegraphics[width=\linewidth]{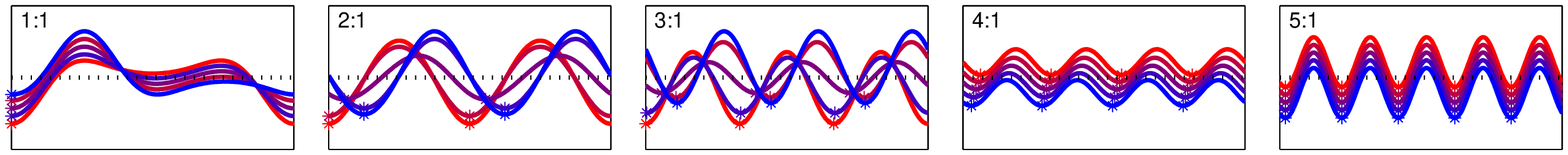}} \vspace{-.25cm}
\caption{\footnotesize Interaction functions for the controls shown in Figure \ref{s6f3_wk}, where the domain and range in each plot is $[0,2\pi]$ and $[-0.05,0.05]$, respectively.  The entrainment ratio is indicated, and the black line denotes the $x$-axis.} \vspace{.25cm} \label{s6f4_winf}
\end{figure}
It follows that to maximize the entrainment range, the phase $\Delta\vphi$ must minimize $V^{\,\uNM}$ in order to maximize the span of the interval $R_\O^{\,\uNM}[v]$, hence  $\Delta\vphi=\vphi_*^{\uNM}$ as in \eqref{magang1}.  In addition, substituting the candidate solution \eqref{sol:candw} into the constraint \eqref{eq:conw1}, we obtain
\begin{align} \label{multw1}
P =\bt{v_r^2} =\frac{1}{4\lambda^2}\bq{2\bt{Y^{\,\uNM}Y^{\,\uNM}}-2\bt{Y^{\,\uNM}(\eta,\vphi^+)Y^{\,\uNM}(\eta,\vphi^-)}} =  \frac{1}{2\lambda^2}\bq{V_0^{\,\uNM}-V_*^{\,\uNM}},
\end{align}
so that solving for the multiplier $\lambda$ yields
\begin{align} \label{multw2}
\lambda = \frac{1}{\sqrt{2P}}\sqrt{V_0^{\,\uNM}-V_*^{\,\uNM}}.
\end{align}
Therefore the waveform of energy $P$ with maximum locking range $R_\O^{\,\uNM}[v]$ for an ensemble of the form \eqref{en_sys_1} is given by
\begin{align} \label{solw2}
 v_r(\eta)=\frac{\sqrt{P}}{\sqrt{2(V_0^{\,\uNM}-V_*^{\,\uNM})}}\bq{Y^{\,\uNM}(\eta,\vphi_*^{\uNM})-Y^{\,\uNM}(\eta,0)}.
\end{align}
Note that although a phase ambiguity exists because we have solved for $\Delta\vphi$, but not for $\vphi^+$ and $\vphi^-$, the initial phase at which the control is applied is unimportant because entrainment is an asymptotic process.  The waveform \eqref{solw2} is actually a special case of \eqref{asoltote} when $\Omega=\half(\w_1+\w_2)$, and the extremal detunings $\Delta\w_2=-\Delta\w_1$ are related to the control energy by $P=2(\Delta\w_2)^2/(V_0^{\,\uNM}-V_*^{\,\uNM})$.  Such controls are shown in {\color{rb3} purple} in Figure \ref{s6f3_wk}.  We may deduce that the control \eqref{solw2} results in the greatest locking range for a fixed control energy, and can be applied at subharmonic forcing frequency $\Omega_f=\frac{N}{M}\half(\w_1+\w_2)$ to entrain the ensemble $\cF$ with minimum control energy. These dual objectives are optimized by the same waveform, and this link  clarifies the relationship between the interaction function and the maximal frequency locking range, which was first observed for harmonic entrainment \cite{harada10,zlotnik11}.

\section{Simulations of Minimum-Energy Entrainment}  \label{sechh}

In this section, we present the results of several numerical simulations that validate the theoretical results that we have derived above for minimum-energy subharmonic entrainment.  Specifically, we compare the theoretical Arnold tongues for the waveforms that were derived from the phase-reduced Hodgkin-Huxley system with computed Arnold tongues for the phase-reduced and full state-space systems. We first apply the phase reduction procedure described in Appendix \ref{cprc} to the equations given in Appendix \ref{aphh}.  Because of the periodicity of $Z$, $v$, and $\Lambda_v^{\uNM}$, all of these functions are conveniently represented using Fourier series, as described in Appendix \ref{apnum}.  These representations are used to synthesize optimal waveforms, which are then applied to simulations to test for entrainment of ordinary differential equation systems for phase models and state-space systems.  Numerical integrations are performed using the $4^{th}$ order Runge-Kutta method.

\begin{figure}[t]
\centerline {\includegraphics[width=\linewidth]{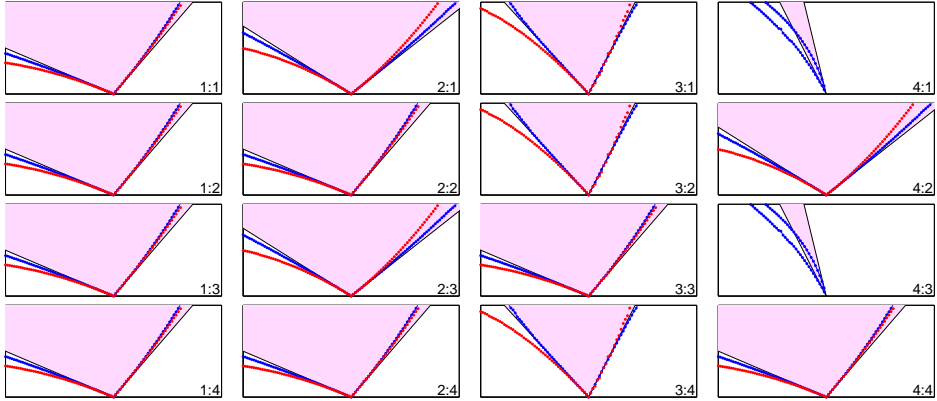}} \vspace{-.25cm}
\caption{\footnotesize Arnold tongues for minimum energy subharmonic entrainment controls $v_m$ for Hodgkin-Huxley neurons, where the target frequency is $\O=0.99\w$.  The domain in each panel is the forcing frequency $\O_f$ on the interval of $90\%$ to $110\%$ of $\frac{N}{M}\w$ where $\w$ is the natural frequency, and the range is $[0,1]$.  The entrainment ratio is indicated at the bottom right.  The shaded region is the theoretical Arnold tongue as determined by Table \ref{s3t1_arnold}.  The actual boundaries of the Arnold tongues are computed as well for entrainment of the phase model, as shown in {\color{blue} blue}, and for the full state-space model, as shown in {\color{red} red}.  For the 4:1 and 4:3 cases, the tongues become too narrow to compute for the state-space model.} \vspace{.25cm}  \label{s8f1_at}
\end{figure}

\begin{figure}[t]
\centerline {\includegraphics[width=1.2\linewidth]{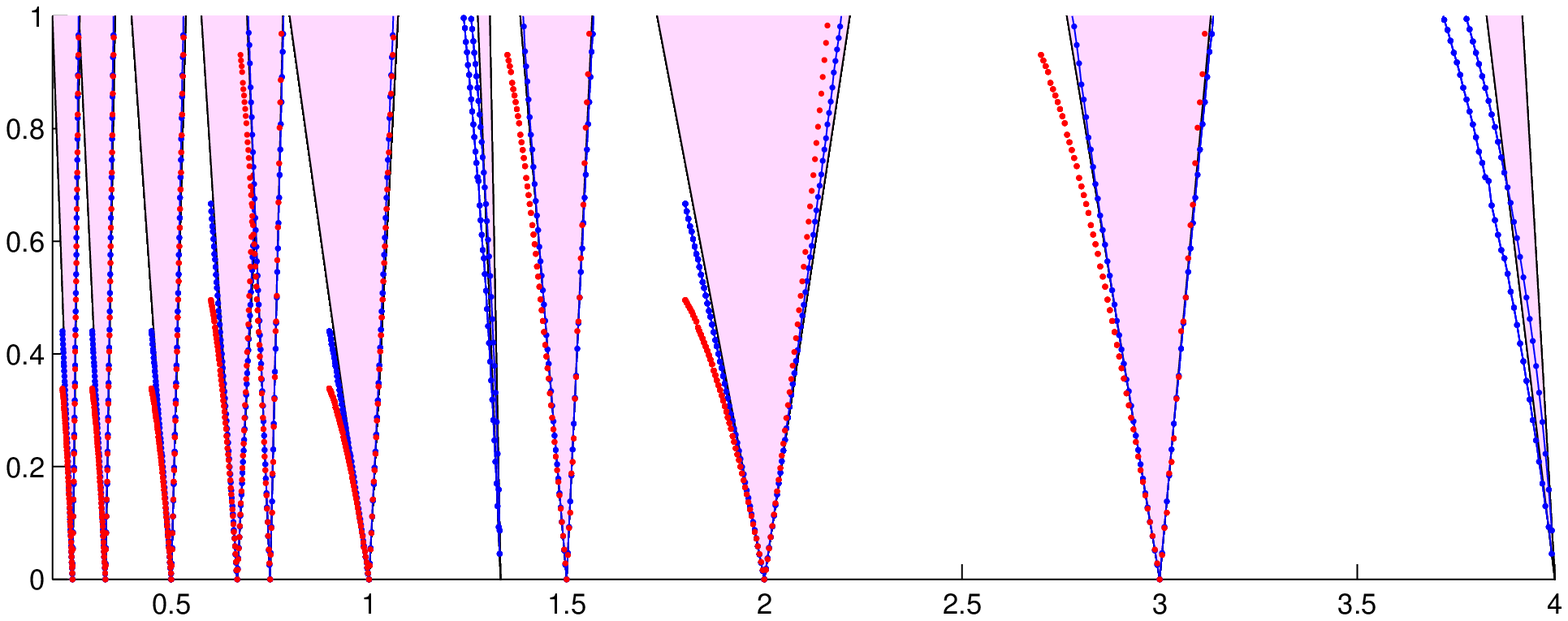}} \vspace{-.25cm}
\caption{\footnotesize Arnold tongues for minimum energy subharmonic entrainment controls $v_m$ for Hodgkin-Huxley neurons, where the target frequency is $\O=0.99\w$.  All of the Arnold tongues in Figure \ref{s8f1_at} are shown together on one plot, where the domain is $\O_f/\w$, i.e., the ratio between the forcing and natural frequencies. } \vspace{.25cm} \label{s8f2_at}
\end{figure}

For the waveform $v_-$ as in \eqref{solnm1}, which is used to entrain a single oscillator, the results are shown in Figures \ref{s8f1_at} and \ref{s8f2_at}.  In this case the natural frequency $\w$ of the oscillator is fixed, and the minimum RMS energy $P_v^{\,\uNM}(\O_f)$ is obtained as a function of the forcing frequency $\O_f$.  The theoretical Arnold tongue is computed by rearranging \eqref{arnold_tg1}, while the actual Arnold tongues for the phase-reduced and state-space system are computed by fixing values of $\O_f$ and using a line search to compute the boundary of the entrainment region.  A bisection search is initialized using guesses of $.9$ and $1.1$ times the theoretical estimate of $P_v^{\,\uNM}(\O_f)$, and is terminated when the upper and lower bound are within $0.01$ times that estimate.  To determine whether a unit energy waveform $\tw{v}$ entrains a phase model \eqref{sys2} for a given pair $(\O_f,P_{\tw{v}}^{\,\uNM})$ of forcing frequency and control energy, the control input $u(t)=P_{\tw{v}}^{\,\uNM}\cdot\tw{v}(\O_f t)$ is applied to the phase model, which is initialized at a fixed point of $\Delta \w + \Lambda_v^{\uNM}(\vphi)$, such as $\vphi_1^*$ in the illustration in Figure \ref{s3f1_inter1}.  The system \eqref{sys2} is integrated numerically, and then the time-series $\psi_k:=\psi(k T_e)$, $k=1,2,\ldots$, where $T_e=2\pi/\O$ is the desired period for the entrained system, is examined to check for convergence to a steady state value.  Convergence of this time-series implies that the forced system has the desired period $T_f$.  In practice, we check whether $\psi_k$ for $k=46,\ldots,50$ remains within an error tolerance of $\ep_1=10^{-1}$.  This approach provides enough time to guarantee that the system has converged to steady-state, in the case that entrainment occurs.  Our experiments have shown that this straightforward approach is sufficient to approximate the minimum RMS energy $P_v^{\,\uNM}(\O_f)$ with error below $1\%$ of the actual value.  We extend the same technique to compute Arnold tongues for the full state-space system \eqref{sys1} by applying to it the same control, and examining the time-series $y_k:=x_1(k T_e)$, $k=1,2,\ldots$, where $x_1(t)$ is the first state variable.  To obtain a reasonably accurate estimate of the boundary, we accept that convergence has occurred when $y_k$ for $k=200,\ldots,250$ remains within an error tolerance of $\ep_2=10^{-2}$.

\begin{figure}[t]
\centerline {\includegraphics[width=\linewidth]{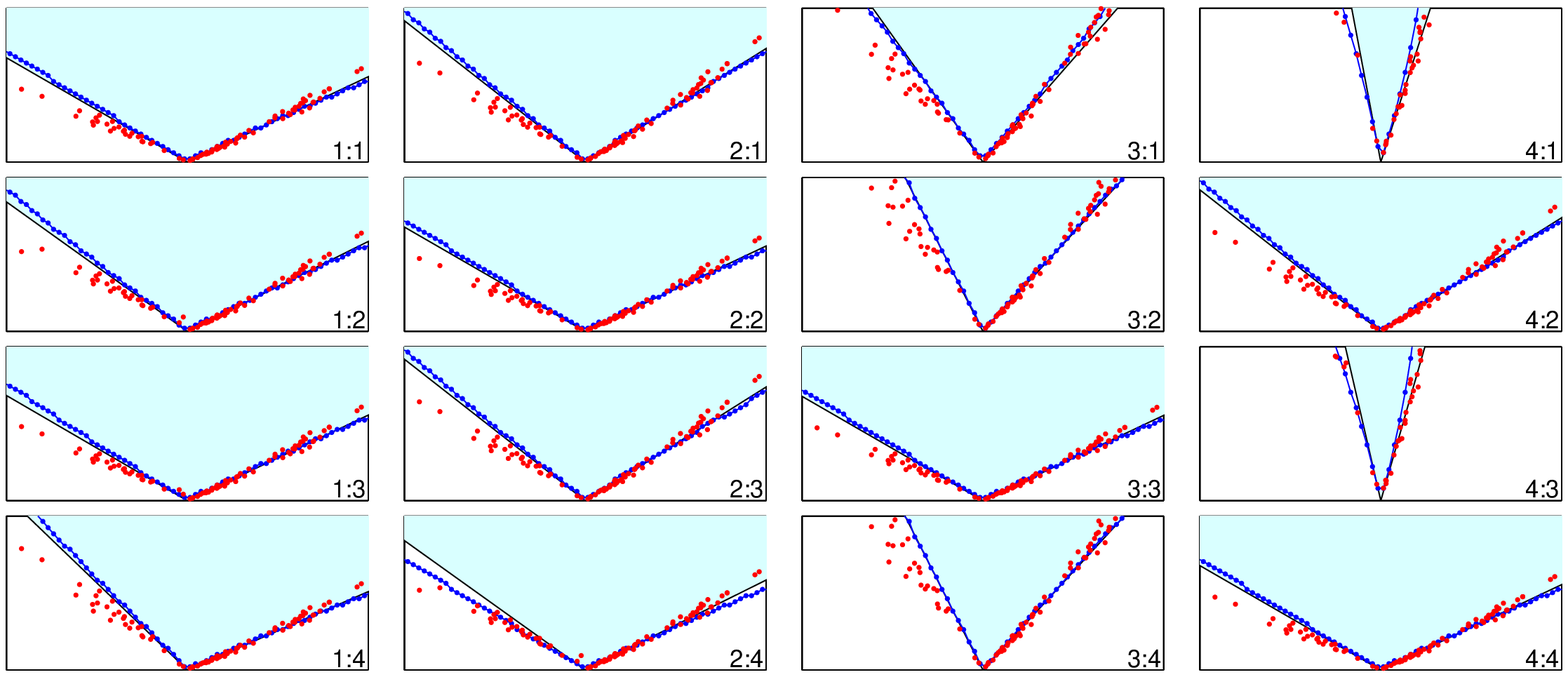}} \vspace{-.25cm}
\caption{\noindent \footnotesize Arnold tongues for minimum-energy subharmonic entrainment ensemble controls $v_e$ for Hodgkin-Huxley neurons, where the target frequency is $\O=0.99\w$.  The domain in each panel is the forcing frequency $\O_f$ on the interval of $90\%$ to $110\%$ of $\frac{N}{M}\w$ where $\w$ is the natural frequency, and the range is $[0,.5]$.  The entrainment ratio is indicated at the bottom right.  The shaded region is the theoretical Arnold tongue as determined by \eqref{ensatng1}.  The computed boundaries of the Arnold tongues are shown in {\color{blue} blue} for entrainment of the phase model, and minimum entrainment energies for the state-space model with parameter values at the corner points of $\cD$ are shown in {\color{red} red}.  } \vspace{.25cm} \label{s8f3_at}
\end{figure}

\begin{figure}[t]
\centerline {\includegraphics[width=1.2\linewidth]{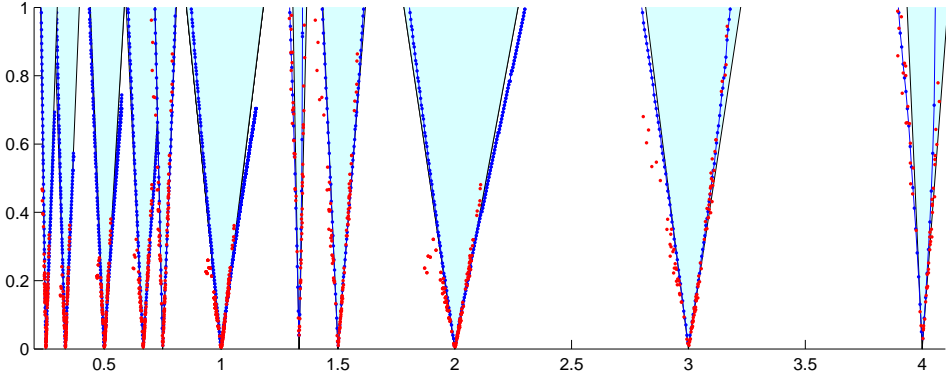}} \vspace{-.25cm}
\caption{\footnotesize Arnold tongues for minimum energy subharmonic entrainment ensemble controls $v_e$ for Hodgkin-Huxley neurons, where the target frequency is $\O=0.99\w$.  All of the Arnold tongues in Figure \ref{s8f3_at} are shown together on one plot, where the domain is $\O_f/\w$, i.e., the ratio between the forcing and natural frequencies. } \vspace{.25cm} \label{s8f4_at}
\end{figure}

Examination of entrainment regions for the waveforms $v_e$ in \eqref{asoltot} that entrain ensembles of oscillators is complicated by the alternative notion of entraining an ensemble, as illustrated in Figure \ref{s6f2_arnold}.  Rather than varying the forcing frequency to compute Arnold tongue of a single oscillator with fixed natural frequency, this notion of an ensemble Arnold tongue requires the forcing frequency to remain fixed while the forcing energy required to entrain oscillators with varying natural frequency is determined.  Recall that in Section \ref{secnmens} we considered a collection of systems $\dot{x}=f(x,u,p)$ where $p\in\cD\subset\bR^d$ is a vector of constant parameters varying on a hypercube $\cD$ containing a nominal parameter vector $q$.  This collection was reduced to the ensemble $\cF$ in \eqref{en_sys_1}.  Thus the ensemble Arnold tongues for the phase-reduced ensemble can be computed by varying $\w$ in $\cF$ and computing the boundaries as described above.  For the ensemble of state-space oscillators, we consider the parameter hypercube $\cD:=\prod_{i=1}^7[.98q_i,1.02q_i]$, where $q=(q_1,\ldots,q_7)$ represents the nominal set of parameters $V_{Na}$, $V_K$, $V_L$, $\ol{g}_{Na}$, $\ol{g}_K$, $\ol{g}_L$, and $c$ of the Hodgkin-Huxley system. Each corner of $\cD$ corresponds to a frequency of oscillation, for which we find the minimum power $P_{\tw{v}}^{\,\uNM}(\w)$ that results in entrainment.  These points are plotted along with the shaded theoretical Arnold tongues in Figures \ref{s8f3_at} and \ref{s8f4_at}, which arise from simulations in which the target frequency is set to $99\%$ of the nominal frequency for the Hodgkin-Huxley system.

Figures \ref{s8f1_at} and \ref{s8f2_at} show close agreement between the phase-locking regions predicted from the theory and the computed boundaries when the forcing frequency $\O_f$ is within several percent of $\frac{N}{M}\w$, and the first order approximation for the phase dynamics given by the model \ref{sys2} is accurate.  For larger frequency detuning, the nonlinear behavior of the oscillation is not captured by the phase model.  Figures \ref{s8f3_at} and \ref{s8f4_at} show that the theoretical and computed ensemble Arnold tongues agree as well.

\section{Simulations of Fast Entrainment}  \label{secfastsim}

In this section, we apply the optimal waveforms for fast subharmonic entrainment given by \eqref{sol:opt} to the Hodgkin-Huxley model in order to verify the performance.  We compare the observed entrainment rates near the asymptotic value of the slow phase with the theoretical value predicted by the gradient of the interaction function, as illustrated in Figure \ref{s5f0_inter1}.  We first note that near the attractive fixed point $\vphi^*$, we can model the dynamics of the averaged slow phase $\vphi$, which evolves according to \eqref{sys5}, using a first order Taylor series approximation about $\vphi^*$.  This takes the form
\begin{align} \label{phas_Tayl}
\dot{\vphi} & = \Delta\w + \Lambda_{v}^{\uNM}(\vphi^*) + \ddx{\vphi}\Lambda_{v}^{\uNM}(\vphi^*) (\vphi-\vphi^*) = \ddx{\vphi}\Lambda_{v}^{\uNM}(\vphi^*) (\vphi-\vphi^*),
\end{align}
where $\Delta\w + \Lambda_{v}^{\uNM}(\vphi^*)\equiv 0$ because $\vphi^*$ is the fixed point that yields $\dot{\vphi}=0$ for \eqref{sys5}.  Setting $\vphi_0 = \vphi-\vphi^*$, the equation \eqref{phas_Tayl} becomes
\begin{align}\label{exp_dec_vphi}
\dot{\vphi}_0&=\ddx{\vphi}\Lambda_{v}^{\uNM}(\vphi^*)\vphi_0
\end{align}
when $\vphi_0$ is near zero.  Recall that the slow phase itself is defined by $\phi(t)=\psi(t)-\O t=\psi(t)-\frac{M}{N}\O_f t$, and follows the dynamics \eqref{sys3}.  Due to the weak forcing assumption, \eqref{sys3} can be approximated near the steady state value $\phi^*$ by \eqref{exp_dec_vphi} where $\vphi_0$ is replaced with $\phi_0=\phi(t)-\phi^*$.  Hence the slow phase $\phi$ decays exponentially to $\phi^*$ according to $\ln|\phi-\phi^*|=c_0+\ddx{\vphi}\Lambda_{v}^{\uNM}(\vphi^*) t$, where $c_0$ is independent of time.  This leads to the following methods for approximating entrainment rates when simulating phase models and state-space systems.

For simulations involving the phase model, we examine the slow phase by simulating the phase model \eqref{sys2} where $u$ is the subharmonic fast entrainment input given by \eqref{sol:opt}.  In particular, we create a time-series $\phi_k:=\phi(k T_e)=\psi(kT_e)-\O k T_e$, $k=1,2,\ldots$, that samples the slow phase system \eqref{sys3}, where $T_e=2\pi/\O$ is the desired period for the entrained system.  The behavior of the phase difference $\phi_k-\phi^*$ in the neighborhood of $\phi^*$ can be closely described by an exponential decay
\begin{align} \label{phas_exp1}
\ln|\phi_k-\phi^*|&=c_0+\kappa_1 k T_e,
\end{align}
where $\kappa_1$ is a negative coefficient that quantifies entrainment rate for the phase model.  Alternatively, \eqref{phas_exp1} leads to the relation
\begin{align} \label{phas_exp2}
\ln|\phi_{k+1}-\phi_k|&=c_1+\kappa_1 k T_e,
\end{align}
where $c_1$ is independent of $k$.

When simulating entrainment of the state-space model, the slow phase $\phi(t)$ must be approximated by locating the peaks of the first state variable $x_1$.  We first form the time-series $z_j:=x_1(t_j)$, $j=1,2,\ldots$, where $t_j$ is the time of the $j^{th}$ peak.  Recall that we define $\psi(t)=0$ (mod $2\pi$) to occur when $x_1$ attains a peak in its cycle, so that $\psi(t_j)=2\pi j$.  We can then define a new slow phase sequence by $\phi_j:=\phi(t_j)=\psi(t_j)-\O t_j$, which yields $t_j=(2\pi j - \phi_j)/\O$, and hence $t_{j+1}-t_j=T_e-(\phi_{j+1}-\phi_j)/\O$, which yields
\begin{align} \label{stat_exp1}
\phi_{j+1}-\phi_j= 2\pi\frac{T_e-(t_{j+1}-t_j)}{T_e}.
\end{align}
Using the slow phase sequence $\{\phi_j\}$ instead of $\{\phi_k\}$ in \eqref{phas_exp2} and applying \eqref{stat_exp1} yields
\begin{align} \label{stat_exp2}
\ln\bl{2\pi\frac{T_e-(t_{j+1}-t_j)}{T_e}}=c_2+\kappa_2 j T_e,
\end{align}
where $c_2$ is independent of $j$ and $\kappa_2$ is a negative coefficient that quantifies the entrainment rate for the state-space model.

\begin{figure}[t]
\centerline {\includegraphics[width=1.2\linewidth]{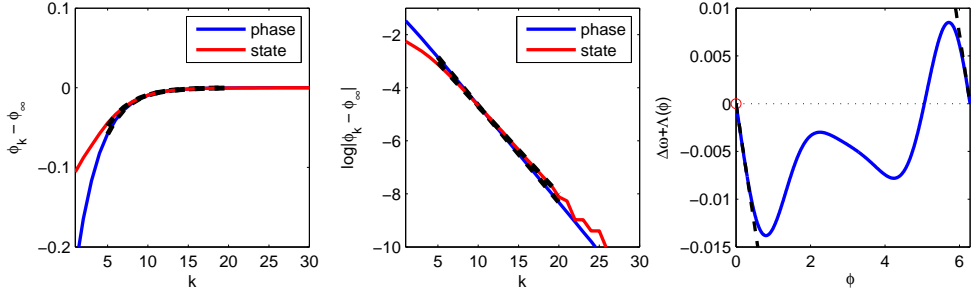}} \vspace{-.25cm}
\caption{\footnotesize Comparing entrainment rates for the Hodgkin-Huxley model.  Left: Phase difference from steady-state value $\psi_k-\psi_\infty$ as a function of period $k$ for harmonic (1:1) fast entrainment simulations of the phase model ({\color{blue} blue}) and state-space model ({\color{red} red}).  Dashed lines show exponential fits.  The target frequency is $\O=1.01\w$.  Center: Linear fits with slopes $\kappa_1$ and $\kappa_2$ of the log of the phase difference from steady state for the phase model and state-space model, according to \eqref{phas_exp1} and \eqref{stat_exp2}, respectively.  The phase converges exponentially to the steady-state.  Right: The dashed line is tangent to $\dot{\vphi}=\Delta\w+\Lambda(\vphi)$ at the attractive phase $\vphi=\vphi^*=0$.  The slope of the line is $\ddx{\vphi}\Lambda^{\uNM}_v(\vphi^*)$, which is the theoretical value of the convergence rate when the oscillator is in the neighborhood of $\vphi^*$.  For the simulations, the system is initialized so that $\phi(0)=\psi^*-0.4$ radians. In this example $\kappa_1=-0.0252$, $\kappa_2=-0.0229$, and $\ddx{\vphi}\Lambda_{v}^{\uNM}(\vphi^*)=-0.0256$.} \vspace{.25cm} \label{s9f1_fer}
\end{figure}

\begin{figure}[t]
\centerline {\includegraphics[width=\linewidth]{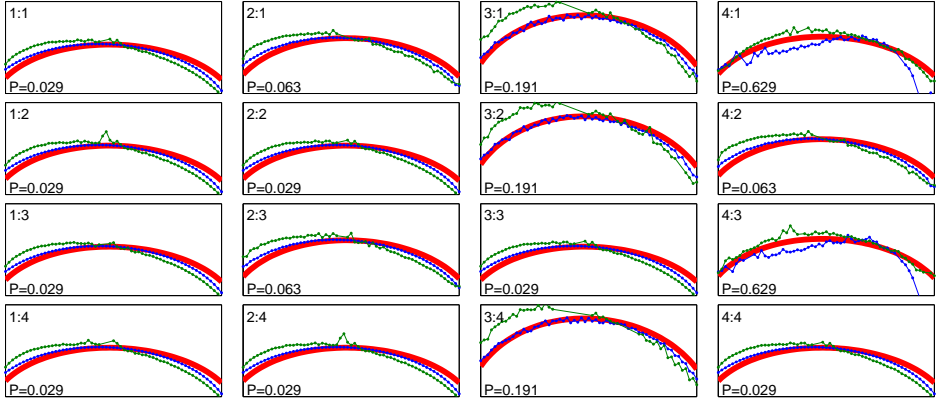}} \vspace{-.25cm}
\caption{\footnotesize Entrainment rates resulting from subharmonic fast entrainment controls $v_f$ in \ref{sol:opt} for the Hodgkin-Huxley neuron model.  The domain in each panel is the forcing frequency $\O_f$ on the interval of $97\%$ to $103\%$ of $\frac{N}{M}\w$ where $\w$ is the natural frequency, and the range is $[0.007,0.050]$.  The entrainment ratio is indicated at the top left, and the control waveform energy, which is adjusted to be slightly greater than the minimum to guarantee entrainment for all the detuning values for each subharmonic ratio, is given at the bottom left.  The red line is the theoretical entrainment rate $\ddx{\vphi}\Lambda_{v}^{\uNM}(\vphi^*)$, and the computed values of $\kappa_1$ and $\kappa_2$ are shown in blue and green, respectively.  } \vspace{.25cm} \label{s9f2_fer}
\end{figure}

The coefficients $\kappa_1$ and $\kappa_2$ are in practice very near to the theoretical entrainment rate $\ddx{\vphi}\Lambda_{v}^{\uNM}(\vphi^*)$, and we expect $\kappa_1$ to be consistently closer to the theoretical value, because the latter is derived from the phase model.  The procedures for obtaining $\kappa_1$ and $\kappa_2$ are illustrated in Figure \ref{s9f1_fer}, which illustrates an example of harmonic (1:1) fast entrainment of the Hodgkin-Huxley system phase model and state-space model where $\O=1.01\w$, and where the theoretical and computed entrainment rates are found to be very similar.  In addition, the same experiment is repeated for subharmonic ($N$:$M$) entrainment and for a range of values of the detuning $\Delta\w$, and the results are given in Figure \ref{s9f2_fer}.  The values are in close agreement, although the computation becomes problematic at higher entrainment ratios.  Observe that the entrainment rate is highest near the center of each panel in Figure \ref{s9f2_fer}, which corresponds to $\Delta\w\equiv 0$.  This is because the frequency of the oscillator does not need to be altered, so that the entrainment rate maximization objective \eqref{prob:1} takes precedence in the problem formulation posed in \eqref{prob:1}-\eqref{eq:con2}.  Conversely, when the detuning is greater, i.e., when $\O$ is farther from $\w$, the optimal theoretical and observed entrainment rate is lower, because the the frequency entrainment design constraint \eqref{eq:con2} influences the problem significantly.

\section{Discussion} \label{secdisc}

The effect of the subharmonic forcing ratio on the entrainment properties of a phase model and given input waveform can be inferred directly from the definition of the interaction function in \eqref{Lambdef}.  The locking range of a control waveform, and hence the Arnold tongue, depends on the properties of the interaction function, as illustrated in Figure \ref{s3f1_inter1} and Table \ref{s3t1_arnold}.  The cases of sinusoidal PRC and sinusoidal forcing are particularly useful.  If the PRC is $Z(\theta)=\sin(\theta)$, then it is evident that $N$:$M$ entrainment is not possible when $N>M$, because the orthogonality of the trigonometric basis functions of the Fourier series would result in $\Lambda_{v}^{\uNM}(\vphi)=0$.  Conversely, if the input is $v(\theta)=\sin(\theta)$, then $N$:$M$ entrainment is not possible when $N<M$ for the same reason.  This leads to the lemma regarding the existence of subharmonic locking regions which is given in Appendix \ref{apnum}.

Furthermore, observe also that as $N$ increases in Figure \ref{s4f1_mpif}, the interaction function quickly narrows to a very small range, which no longer contains the origin, so that the Arnold tongue will skew to one side.  This is illustrated in Case C of Figure \ref{s3f2_arnold}, and is observed in practice for the $4$:$1$ and $4$:$3$ cases in Figure \ref{s8f1_at}.  This is due to the rapid decrease of energy in successive terms of the Fourier series for the Hodgkin-Huxley PRC.  In fact, $N$:$M$ entrainment can be established only if there is significant energy in the $N^{th}$ term of this series. Because the coefficients of the Fourier series for the Hodgkin-Huxley PRC are nearly negligible beyond the fourth order, the Arnold tongues also become extremely thin, and so that subharmonic entrainment with $N=5$ (and $M<5$) cannot be established in practice for this system.

Throughout this paper, we focus on deriving waveforms which are optimal in the case of the weak forcing assumption, as described in Section \ref{secent}.  Many approximations are made in the process of phase reduction and averaging, so that the controls presented above are optimal only in an approximate sense, as the input energy $\bt{v^2}$ approaches zero.  An analysis of the accuracy and divergence from optimality of the produced controls, taking into account control amplitude, accuracy of phase reduction, and effects of averaging is a challenging problem that is left for future work.  It is important to emphasize that the techniques presented here provide a straightforward way to compute near optimal controls numerically, and have been applied successfully in an experimental setting \cite{zlotnik13prl}.  Furthermore, from the simulation shown in Figure \ref{s8f4_at}, we see that the optimal waveforms obtained using the phase modeling technique produce a similar result to the theory when applied to both the phase model and original model.  This strongly supports the hypothesis that optimal entrainment controls derived using a phase model are very near optimal for the original system, provided the oscillator remains within a neighborhood of its limit cycle.

The results on subharmonic entrainment of oscillator ensembles in Section \ref{secnmens} can be interpreted as a means of shaping the Arnold tongue characterizing the entrainment of an oscillatory ensemble.  By adjusting the forcing waveform while keeping the forcing frequency $\O_f$ fixed, it is possible to significantly alter the frequency range $(w_1,w_2)$ of the collection of oscillators $\cF$ subjected to subharmonic phase-locking.  The main focus here is on accomplishing such manipulation in an optimal manner.  In most cases, entrainment of a given ensemble $\cF$ can actually be achieved using a biased sinusoid of the form $u(t)=b_0+b_1\sin(\O_f t)$ with appropriate constants $b_0$, $b_1$, and $\O_f$.  However, the derived waveforms accomplish this design goal using significantly less energy.  Furthermore, the analysis of Arnold tongues for minimum-energy waveforms provides a framework for studying the possibilities and limitations of engineering entrainment of rhythmic systems on multiple time scales, for instance in an interacting network.  Such analysis may also shed light on the evolved optimal periodic activity of complex multi-scale biological systems.  For example, experimentally measured subharmonic entrainment regions were approximated by injecting single Aplysia motoneurons with sinusoidal inputs of varying frequency and amplitude \cite{hunter03}, resulting in plots similar to Figures \ref{s8f3_at} and \ref{s8f4_at}.

An indirect implication of this work stems from the importance of the interaction function between the PRC of one or more oscillators and the common control input for the phase-locking properties of the ensemble.  In this paper we have focused entirely on the frequency locking aspect of entrainment, and did not consider the fixed point of the average slow phase $\vphi$ to which oscillators converge, as long as the frequency control objective is satisfied.  However, as one can see in Figure \ref{s3f1_inter1}, the form of the interaction function determines the asymptotic phases of entrained oscillators, which is of particular interest when manipulating the synchronization of multiple rhythmic units.  The techniques presented here can therefore be extended to engineer synchronization in collections of oscillators using weak forcing without feedback information, which is of compelling interest in electrochemistry \cite{kiss07}, neuroscience \cite{uhlhaas06}, and circadian biology \cite{gooley10}.  The impact may be greatest on the ability to manipulate collections of biological circadian and neural systems, for which the entrained phases may need to be design in a nonuniform manner, by using common inputs.  Novel paradigms for designing synchronization patterns in rhythmic biological and electrochemical systems will appear in our future work.

%A study of entrainment in strongly forced systems is a compelling extension \cite{kurebayashi13}.

\section*{Conclusions} \label{secconc}

We have developed a methodology for designing optimal waveforms for subharmonic entrainment of oscillator ensembles to a desired frequency using weak periodic forcing.  Our approach is based on  phase model reduction, formal averaging theory, and the calculus of variations.  Diverse objectives such as minimizing input energy, maximizing the rate of entrainment, and designing the entrainment frequency and control power are considered.  In addition, the entrainment of large ensembles of oscillators with uncertain parameters is explored.

In order to characterize the phenomenon of subharmonic entrainment, we also derive an approximation of the locking region in energy-frequency space for a periodically forced oscillator, called an Arnold tongue, or for an oscillator ensemble, which we refer to as an ensemble Arnold tongue.  The entrainment of phase-reduced Hodgkin-Huxley neurons is considered as an example problem throughout the paper, and boundaries of Arnold tongues are computed for various subharmonic entrainment ratios and controls to compare to the theoretical regions.  Detailed descriptions and illustrations are provided to connect the behavior of entrained systems to the corresponding Arnold tongues and interaction functions.  Simulations to compute actual Arnold tongues are described and carried out for minimum energy subharmonic entrainment controls for single oscillators and oscillator ensembles, and the measurement of entrainment rates of simulations is described and carried out to examine the performance of fast entrainment waveforms.  In all cases, the computational results closely correspond to what is predicted by the derived theory.

This work provides a comprehensive study of subharmonic entrainment of weakly forced nonlinear oscillators, as well as a practical technique for control synthesis.  The methods presented here may also be used for the analysis of synchronization in interacting rhythmic systems across different time-scales.  The approach described is of direct interest to researchers in chemistry and biology, and particularly neuroscience.

\appendix

\section{Interaction functions for subharmonic entrainment} \label{apnum}

Because the PRC $Z(\theta)$, input waveform $v(\theta)$, and interaction function $\Lambda_v^{\uNM}(\vphi)$ are all $2\pi$-periodic, they are most conveniently represented using Fourier series, and interaction functions can easily be computed by inspecting the equation \ref{Lambdef}.  Let us denote the Fourier series for $Z$ and $v$ by
\begin{align}
\dS Z(\theta) &= \dS \half a_0 + \sum_{n=1}^\infty a_n\cos(n\theta) + \sum_{n=1}^\infty b_n\sin(n\theta),\\
\dS v(\theta) &= \dS \half c_0 + \sum_{n=1}^\infty c_n\cos(n\theta) + \sum_{n=1}^\infty d_n\sin(n\theta).
\end{align}
By applying trigonometric angle sum identities and the orthogonality of the Fourier basis, we can obtain $\Lambda_v^{\uNM}(\vphi)$ by first computing
\begin{align} \label{applamb}
\Lambda_v^{\uNM}(M\vphi) &= \dS \frac{a_0c_0}{4}  + \half \sum_{j=1}^\infty [a_{M\!j\,}c_{N\!j\,}+b_{M\!j\,}d_{N\!j\,}]\cos(j\vphi) \nonumber \\ & \qquad +\half \sum_{j=1}^\infty [b_{M\!j\,}c_{N\!j\,}-a_{M\!j\,}d_{N\!j\,}]\sin(j\vphi),
\end{align}
then making the appropriate re-scaling.  The integers $N$ and $M$ must be coprime.  The equation \eqref{applamb} leads to the following lemma:

\noindent \tb{Lemma 1:} \tb{Condition for existence of subharmonic entrainment.} Given a phase model \eqref{sys2} and an input waveform $v(\theta)$, subharmonic ($N$:$M$) entrainment to a target frequency $\O\neq \w$ using a forcing frequency $\O_f=\frac{N}{M}\O$ is possible if and only if $a_{M\!j\,}c_{N\!j\,}+b_{M\!j\,}d_{N\!j\,}\neq 0$ or $b_{M\!j\,}c_{N\!j\,}-a_{M\!j\,}d_{N\!j\,}\neq 0$ for at least one $j\in\bN$.

\noindent \tb{Proof:} If $a_{M\!j\,}c_{N\!j\,}+b_{M\!j\,}d_{N\!j\,}= 0$ and $b_{M\!j\,}c_{N\!j\,}-a_{M\!j\,}d_{N\!j\,}= 0$ for all $j\in\bN$, then $\Lambda_v^{\uNM}(\vphi)\equiv0$ by \eqref{applamb}.  Therefore if $\Delta\w=\w-\O\neq 0$, then $\Lambda_v^{\uNM}(\vphi)+\Delta\w=0$ has no solution and \eqref{sys5} has no fixed point.  Therefore entrainment cannot occur.

Conversely, suppose that $N$:$M$ entrainment is possible for $\Delta\w\neq 0$.  It follows that $\Lambda_v^{\uNM}(\vphi)+\Delta\w=0$ must have a solution, hence $\Lambda_v^{\uNM}(\vphi)$ is not identically zero, and therefore $a_{M\!j\,}c_{N\!j\,}+b_{M\!j\,}d_{N\!j\,}\neq 0$ or $b_{M\!j\,}c_{N\!j\,}-a_{M\!j\,}d_{N\!j\,}\neq 0$ for at least one $j\in\bN$. \boX

\noindent \tb{Lemma 2:} \tb{Continuity of interaction function.} Interaction function $\Lambda_v^{\uNM}$ is continuous when $v$ is bounded.

\noindent \tb{Proof:} Suppose $v$ is bounded but $\Lambda_v^{\uNM}$ in \eqref{Lambdef} is discontinuous at $\phi\in[0,2\pi)$.  Then $\eX M>0$ such that $|v(\theta)|<M$ $\fA \theta\in[0,2\pi)$, and $\eX \sigma>0$ such that $|\Lambda_v^{\uNM}(\phi+\frac{\delta}{2})-\Lambda_v^{\uNM}(\phi-\frac{\delta}{2})|>\sigma$ $\fA \delta\in(0,d)$ for some $d>0$. Because $Z$ is continuous, it follows that $\fA \ep>0$, $\eX \delta>0$ such that $|Z(\theta+\delta)-Z(\theta)|<\epsilon$.  Then for arbitrary fixed $\eP>0$,
\begin{align*}
0& <\sigma<|\Lambda_v^{\uNM}(\phi+\tfrac{\delta}{2})-\Lambda_v^{\uNM}(\phi-\tfrac{\delta}{2})| \\
& = \bl{\frac{1}{2\pi}\int_0^{2\pi} [Z(M\theta+\phi+\tfrac{\delta}{2})-Z(M\theta+\phi-\tfrac{\delta}{2})]v(N\theta)\rd\theta} \\
& \leq \frac{1}{2\pi}\int_0^{2\pi} |Z(M\theta+\phi+\tfrac{\delta}{2})-Z(M\theta+\phi-\tfrac{\delta}{2})|\cdot|v(N\theta)|\rd\theta \\ &\leq \frac{1}{2\pi}\int_0^{2\pi} \eP M \rd\theta = \eP M.
\end{align*}
Therefore $\eP M>\sigma$, and choosing $\eP=\sigma/(2M)$ yields a contradiction. \boX

\section{Hodgkin-Huxley Model} \label{aphh}

The Hodgkin-Huxley model describes the propagation of action potentials in neurons, specifically the squid giant axon, and is used as a canonical example of neural oscillator dynamics. The equations are
\begin{eqnarray*} \label{hheq}
\hskip-2pt\begin{array}{c}\begin{array}{rcl}
c\dot{V}&=&I_b+I(t)-\ol{g}_{Na}h(V-V_{Na})m^3-\ol{g}_K(V-V_k)n^4-\ol{g}_L(V-V_L)\\
\dot{m}&=&a_m(V)(1-m)-b_m(V)m,\\
\dot{h}&=&a_h(V)(1-h)-b_h(V)h,\\
\dot{n}&=&a_n(V)(1-n)-b_n(V)n,
\end{array}\\\\\begin{array}{rcl}
a_m(V)&=&0.1(V+40)/(1-\exp(-(V+40)/10)),\\ b_m(V)&=&4\exp(-(V+65)/18),\\
a_h(V)&=&0.07\exp(-(V+65)/20),\\ b_h(V)&=&1/(1+\exp(-(V+35)/10)),\\
a_n(V)&=&0.01(V+55)/(1-\exp(-(V+55)/10)),\\ b_n(V)&=&0.125\exp(-(V+65)/80).\\
\end{array}\end{array}
\end{eqnarray*}
The variable $V$ is the voltage across the axon membrane, $m$, $h$, and $n$ are the ion gating variables, $I_b$ is a baseline current that induces the oscillation, and $I(t)$ is the control input.  The units of $V$ are millivolts and the units of time are milliseconds. Nominal parameters are $V_{Na} =50 \text{ mV}$, $V_K=-77 \text{ mV}$, $V_L=-54.4 \text{ mV}$, $\ol{g}_{Na}=120 \text{ mS/cm}^2$, $\ol{g}_K = 36 \text{ mS/cm}^2$, $\ol{g}_L=0.3 \text{ mS/cm}^2$, $I_b=10 \,\,\mu\text{A/cm}^2$, and $c=1 \,\,\mu\text{F/cm}^2$, for which the period of oscillation is $T=14.63842\pm10^{-5}$ ms.

\section{Computation of Phase Response Curves} \label{cprc}

In this appendix we present a basic summary of the technique for phase coordinate transformation, derived directly from the method of Malkin \cite{malkin49}.  This derivation leads to a straightforward method for numerical computation of phase response curves, which is implemented here to automatically compute phase models from the Hodgkin-Huxley equations for numerous parameter sets in order to produce Figures \ref{s8f3_at} and \ref{s8f4_at} above.  A collection of theorems regarding existence, accuracy, and validity of phase-reduced models has been produced \cite{efimov11}.

Consider a smooth ODE system $\Dx=f(x,u)$, where $x(t)\in\bR^n$ is the state and $u(t)\in\bR$ is a control, such that $\Dx=f(x,0)$ has an attractive, non-constant, $T$-periodic limit cycle $\gamma(t)=\gamma(t+T)\in\Gamma\subset\bR^n$ evolving on the periodic orbit $\Gamma=\set{y\in\bR^n}{y=\gamma(t) \text{ for } 0\leq t< T}$.  A bijection can be defined between $\Gamma$ and the circle $S^1$, which is homeomorphic to the interval $[0,2\pi)$, hence any point $x\in \Gamma$ can be associated with a scalar phase $\phi\in[0,2\pi)$ by a map $\alpha:\Gamma\to[0,2\pi)$ with action $\alpha(x)=\phi$.   We choose $\alpha$ such that the phase is proportional to time on the the limit cycle, i.e., $\alpha^{-1}(\phi)=\gamma(\phi/\w)$, where $\w=2\pi/T$ is the natural oscillator frequency, and so $\gamma(0)=\alpha^{-1}(0)=\gamma(2\pi)$.   Denote by $x(t,x_0,u)$ a trajectory satisfying $\Dx=f(x,u)$ for a control function $u:[0,t]\to\bR$ and $x(0)=x_0$.  Then $\gamma(t)=x(t,\alpha^{-1}(0),0)$, so that if $x(0)=x_0\in\Gamma$ then $x(t,x_0,0)=\gamma(t+\phi_0/\w)$, where $\phi_0=\alpha(x_0)$.  We can define a phase variable $\phi:[0,\infty)\to[0,2\pi)$ for trajectories $x(t,x_0,0)$, $x_0\in\Gamma$ by $\phi(t)=\alpha(x(t,x_0,0))=\alpha(\gamma(t+\phi_0/\w))$.  Because $\gamma(t)$ is periodic, then $\phi(t)$ is periodic, and our choice of $\alpha$ results in an affine system $\phi(t)=\w t+\phi_0$, so that $\dot{\phi}(t)=\w$.  For any $x_0\in\Gamma$, we define $\gamma(0)=x_0$, so that $\phi_0=\alpha(x_0)=\alpha(\gamma(0))=0$.

Denote by $\cA=\set{y\in\bR^n}{\lim_{t\to\infty}x(t,y,0)\in\Gamma}\subset\bR^n$ the set attracted by the periodic orbit $\Gamma$, so if $x_0\in\cA$ then $x(t,x_0,0)\in\cA$ for $t\geq 0$.  This allows us to extend the notion of phase mapping to any solution $x(t,x_0,0)$ for $x_0\in\cA$, by defining an asymptotic phase $\theta_0\in[0,2\pi)$ such that $\lim_{t\to\infty}\norms{x(t,x_0,0)-\gamma(t+\theta_0/\w)}=0$.  We can define an asymptotic phase map $\upsilon:\cA\to[0,2\pi)$ that maps the point $x_0\in \cA$ to the corresponding phase $\theta_0$, i.e., $\theta_0=\upsilon(x_0)$.  In the case that $x_0\in\Gamma$, then $\norms{x(t,x_0,0)-\gamma(t+\phi_0/\w)}=0$, so that $\upsilon(x_0)=\theta_0=\phi_0=\alpha(x_0)$.  The asymptotic phase variable $\theta:[0,\infty)\to[0,2\pi)$ is a mapping $\theta(t)=\upsilon(x(t,x_0,0))$ for $t\geq 0$, which is defined for $x_0\in\cA$.  Therefore if $x_1,x_2\in\cA$ satisfy $\upsilon(x_1)=\upsilon(x_2)$, then $\upsilon(x(t,x_1,0))=\upsilon(x(t,x_2,0))$ for all $t\geq0$.  We define an equivalence class on $\cA$ by $x_1\sim x_2$ if $\upsilon(x_1)=\upsilon(x_2)$, and denote class elements in the quotient space by $[x_0]\in\cA/\sim$, so that $\upsilon([x_0])=\theta_0$.  We call the class element $[x_0]$ the isochron corresponding to the phase $\theta_0$.  Let $x_1\in\cA$ satisfy $x_1\in[x_0]$ for $x_0\in\Gamma$, so that $\upsilon(x_1)=\upsilon(x_0)$.  Then  $\theta(t)=\upsilon(x(t,x_1,0))=\upsilon(x(t,x_0,0))=\alpha(x(t,x_0,0))=\phi(t)=\w t+\alpha(x_0)=\w t+\upsilon(x_0)=\w t+\upsilon(x_1)=\w t + \theta_0$.  Therefore the asymptotic phase along $\gamma(t)$, when $u\equiv 0$, satisfies $\theta(t)=\w t+ \theta_0$ and $\dot{\theta}(t)=\w$.  For any $x_1\in\cA$, there exists $x_0\in \Gamma$ such that $x_1\sim x_0$.  Therefore we can define $\gamma(0)=x_0$ where $x_1\in[x_0]\in\cA/\sim$, so that $\theta_0=\upsilon(x_1)=\upsilon(x_0)=\alpha(x_0)=\alpha(\gamma(0))=0$.  The following diagram displays the relevant mappings and spaces.
$$\xymatrix{t\in[0,T)
\ar[d]_{\gamma}\ar[rrr]^{(\alpha\circ\gamma)} &&& \theta(t)\in[0,2\pi)
\\ \gamma(t)\in\Gamma \ar[rrru]_\alpha &&&
x(t)\in\cA\ar[u]^{\upsilon}\ar[lll]^{(\alpha^{-1}\circ\upsilon)}}
$$
We can extend the notion of asymptotic phase to the case when $u(t)\neq 0$, provided that $x(t,x_0,u)\in \cA$ for $t\geq 0$.  In this case we define a new asymptotic phase map $\nu:\cA\to[0,2\pi)$ that acts according to $\theta(t)=\nu(x(t,x_0,u))=\upsilon(x(t,x(t,x_0,u),0))$, so that $\theta(t)$ at a time $t\geq 0$ evaluates the asymptotic phase of the point $x(t,x(t,x_0,u),0)$.
In other words, it is the asymptotic phase of $x(t,x_0,z)$ where $z(s)=u(s)$ for $s\in[0,t)$ and $z(s)=0$ for $s\geq t$.

It is possible to use a linearization of the system $\Dx=f(x,u)$ about its limit cycle $\gamma(t)$ to obtain the ODE for the asymptotic phase variable $\theta(t)$ given infinitesimal inputs $u(t)$ such that the solution $x(t,x_0,u)$ remains within a neighborhood of the periodic orbit $\Gamma$.  Define the perturbation variable $\Delta x(t)=x(t)-\gamma(t)$, so that the linearization about $\gamma$ is
$ \Delta \dot{x}(t)=A(t)\Delta \dot{x}(t) + b(t) u$ where
$$
A(t)=\evalat{\ppx{x}f(x,0)}{x=\gamma(t)}{} \anD b(t)=\evalat{\ppx{u}f(\gamma(t),u)}{u=0}{}
$$
Note that $A(t)$ and $b(t)$ are $T$-periodic because they depend on $\gamma(t)$, hence we can apply Floquet theory to the linearized system \cite{perko90}.  The fundamental matrix $\Phi(t)$ satisfies $\dot{\Phi}(t)=A(t)\Phi(t)$ with $\Phi(0)=I$, and its adjoint  $\dot{\Psi}(t)=-A^{\dagger}(t)\Psi(t)$ with $\Psi(0)=I$, where $\dagger$ denotes the Hermitian transpose.  Recall that $\Psi^{\dagger}(t)\Phi(t)=\Psi^{\dagger}(0)\Phi(0)\equiv I$, and that $y(t) = \Phi(t)y(0)$ if $\dot{y}(t)=A(t)y(t)$.  Recall also Floquet's theorem \cite{kelley04}, which states that if $A(t)$ is a continuous, $T$-periodic matrix, then for all $t\in\bR$ any fundamental matrix solution $\Phi$ for $\Dx=A(t)x$ can be written in the form $\Phi(t)=Q(t)e^{Bt}$ where $Q(t)$ is a nonsingular, differentiable, $T$-periodic matrix and $B$ is a constant matrix.  Furthermore, if $\Phi(0)=I$ then $Q(0)=I$.

Define the monodromy matrix $M(t)=\Phi(t+T)\Psi^{\dagger}(t)$, which is the linearized return map of the dynamical system.  By Floquet's theorem, there exists a matrix $B$ such that $\Phi(t)=Q(t)e^{Bt}$ where $Q(t)$ is $T$-periodic, hence
\begin{align}
M(t) & = \Phi(t+T)\Psi^{\dagger}(t) =  \Phi(t+T)\Phi^{-1}(t) \nonumber \\ &=  Q(t+T)e^{B(t+T)}e^{-Bt}Q^{-1}(t) =  Q(t)e^{BT}Q^{-1}(t).
\end{align}
Therefore $M(t)$ is $T$-periodic and isospectral, because $M(t)$ is similar to the constant matrix $M(0)=e^{BT}$.    The eigenvalues $\lambda_i$ of $M(t)$ are the Floquet multipliers of the linearization, each of which corresponds to an eigenvalue $\rho_i$ of $B$, called characteristic exponents, where $\lambda_i=e^{\rho_iT}$.  One of the Floquet multipliers of $\gamma(t)$ is always 1, and by the stable manifold theorem the other $n-1$ multipliers are less than 1 \cite{perko90}.

In the case $u(t)=0$, the limit cycle satisfies $\dot{\gamma}(t)=f(\gamma(t),0)$, so that $\ddot{\gamma}(t)=A(t)\dot{\gamma}(t)$ and hence $\dot{\gamma}(t)=\Phi(t)\dot{\gamma}(0)$.  In particular, $ \dot{\gamma}(t)= \dot{\gamma}(t+T)=\Phi(t+T)\dot{\gamma}(0) =\Phi(t+T)\Phi^{-1}(t)\dot{\gamma}(t) =M(t)\dot{\gamma}(t) $.  It follows that $\dot{\gamma}(t)=f(\gamma(t),0)$ is the unique eigenvector of $M(t)$ corresponding to the Floquet multiplier $\lambda=1$, which has algebraic multiplicity of 1.  Let $m(t)\in\bR^n$ be the unique eigenvector of $M^{\dagger}(t)$ corresponding the the Floquet multiplier $\lambda=1$, so that $M^{\dagger}(t)m(t)=m(t)$, and scaled such that $m^{\dagger}(t)\dot{\gamma}(t)=\w$. It follows that
$m^{\dagger}(t)f(\gamma(t),0)=m^{\dagger}(t)\dot{\gamma}(t)=\w$.

Recall that we defined a mapping $\nu:\cA\to[0,2\pi)$ that acts by $\theta(t)=\nu(x(t,x_0,u))$, and which if chosen properly results in a phase model $\dot{\theta}(t)=\w$ when $u=0$.  Then for $x_0\in \Gamma$,
\begin{eqnarray*}
m^{\dagger}(t)f(\gamma(t),0) & = & \w=\dot{\theta}=\ddx{t}\nu(x(t,x_0,0)) \\
& = & \evalat{\ppx{x}\nu(x)}{x=\gamma(t)}{}\cdot \dot{\gamma}(t) = \evalat{\ppx{x}\nu(x)}{x=\gamma(t)}{}\cdot f(\gamma(t),0)
\end{eqnarray*}
We may therefore infer that $m^{\dagger}(t)=\ppx{x}\evalatb{\nu(x)}{x=\gamma(t)}{}$, and deduce that $m(t)$ is $T$-periodic because $\gamma(t)$ is.  When $u(t)\neq 0$, the linearized trajectory satisfies
\begin{eqnarray*}
\dot{\theta}(t) & = & \evalat{\ppx{x}\nu(x)}{x=\gamma(t)+\Delta x(t)}{}\cdot \dot{x}(t) \\  & =  & \evalat{\ppx{x}\nu(x)}{x=\gamma(t)+\Delta x(t)}{}\cdot (f(\gamma(t),0)+A(t)\Delta x(t) + b(t) u)
\end{eqnarray*}
To obtain the infinitesimal PRC we set $\Delta x(t)=0$, so that $\dot{\theta}(t)=\w+m^{\dagger}(t)b(t)u$.   Now when $x(t)=\gamma(t)$ we have shown that $\theta(t)=\w t$, so we use $t= \theta(t)/\w$.  It follows that we can write $m^{\dagger}(t)b(t)=m^{\dagger}(\theta/\w)b(\theta/\w)$.
We then obtain the system $\dot{\theta}(t)=\w+Z(\theta)u$ where $Z(\theta)=m^{\dagger}(\theta/\w)b(\theta/\w)$ is the PRC.

We see that $m(t)=M^{\dagger}(t)m(t)= (\Phi(t+T)\Phi^{-1}(t))^{\dagger}m(t) = \Psi(t)\Psi^{-1}(t+T)m(t) = \Psi(t)m(-T)=\Psi(t)m(0)$, hence $m(t)$ satisfies the adjoint equation $\dot{m}(t)=-A^{\dagger}(t)m(t)$.  This leads to the following method for computing $Z(\theta)$:

\vspace{.25cm}
\noindent \tb{Algorithm 1:} \tb{Adjoint method.}
\begin{enumerate}
\item Choose point $x_0\in\Gamma$ and compute $\gamma(t)=x(t,x_0,0)$ by integrating $\Dx(t)=f(x(t),0)$ with $x(0)=x_0$.
\item Compute $\Phi(T)$ by integrating $\DPhi(t)=A(t)\Phi(t)$ with $\Phi(0)=I$.
\item Compute $M=M(0)=\Phi(T)$ and the eigenvector $\mu$ of $M^T$, and set $m_0=\w(\mu^Tf(x_0,0))^{-1}\cdot \mu$.
\item Integrate $\dot{m}(t)=-A^T(t)m(t)$ with $m(0)=m_0$.
\item The PRC is given by $Z(\theta)=m^T(\theta/\w)b(\theta/\w)$.
\end{enumerate}
\vspace{.25cm}

The above algorithm has issues with computational stability, because integration of the adjoint system is numerically unstable if $\Phi(t)$ is poorly conditioned.  Rather than using the adjoint method, we can employ the projection method which, although computationally costlier, does not have issues of numerical instability and thus does not require optimization using polynomial approximation schemes.

\vspace{.25cm}
\noindent \tb{Algorithm 2:} \tb{Projection Method.}
\begin{enumerate}
\item Compute the limit cycle $\gamma(t)$, the period $T$, and the natural frequency $\w$.
\item For each $\theta\in[0,2\pi)$, let $x_\theta=\alpha^{-1}(\theta)=\gamma(\theta/\w)$, and compute $x(t,x_\theta,0)$ using $\Dx(t)=f(x(t),0)$ for $x(0)=x_\theta$.
\item Let $A_\theta(t)$ be the linearization of $f(x(t),0)$ about $\gamma_\theta(t)=\gamma(t+\theta/\w)$, so that $\gamma_\theta(0)=x_\theta$.  Compute $\Phi_\theta(T)$ where $\DPhi_\theta(t)=A_\theta(t)\Phi_\theta(t)$ with $\Phi_\theta(0)=I$.
\item Set $M_\theta=\Phi_\theta(T)$, and compute the eigenvector $\mu_\theta$ of $M_\theta^{\dagger}$.
\item $m_\theta=\w (\mu_\theta^{\dagger}f(x_\theta,0))^{-1}\cdot \mu_\theta$, and $Z(\theta)=m_\theta^{\dagger}b(\theta/\w)$.
\end{enumerate}
\vspace{.25cm}

Algorithm 2 eliminates the need to integrate the adjoint equations, wherein numerical instability occurs when the adjoint of the linearization is ill-conditioned.  However, the forward equations must be solved for a full cycle to evaluate $Z(\theta)$ at every desired $\theta$ value.  For further details and useful software packages, we refer the reader to several sources  \cite{ermentrout96,govaerts06,ermentrout02}.

%\section{Hindmarsh-Rose Model} \label{aphr}
%
%The Hindmarsh-Rose model is a generic model of bursting neuron dynamics that is widely used in neuroscience \cite{hindmarsh84}.  A slow variable is augmented to the Fitzhugh-Nagumo planar model, which is a reduction of the Hodgkin-Huxley model \cite{fitzhugh61}.  The equations are
%\begin{eqnarray*} \label{hreq}
%\dot{V}&=&I_b+I(t)+n-aV^3+bV^2-h\\
%\dot{n}&=&c-dV^2-n\\
%\dot{h}&=&r(\sigma(V-V_0)-h)
%\end{eqnarray*}
%The non-dimensional variable $V$ is the voltage across the axon membrane, while $n$ and $h$ are fast and slow ion gating variables, respectively, $I_b$ is a baseline current that induces the oscillation, and $I(t)$ is the control input.  The units of $V$ are millivolts and the units of time are milliseconds. Nominal dimensionless parameters are $a=1 $, $b=3$, $c=1$, $d=5$, $r=0.001$, $\sigma=4$, $V_0=-1.6$, and $I_b=2$, for which a full period of oscillation is $T=430.7757$, corresponding to angular frequency $\w \approx 0.0145857$.

\section{Acknowledgements}

This work was supported by the National Science Foundation under the award 1301148.  We thank Istv\'an Z. Kiss for his insight on entrainment phenomena.

\bibliographystyle{unsrt}
\bibliography{nm_bib}

\end{document}

%% file: avz_defs.tex
\usepackage{amsmath, amssymb, amsxtra, mathrsfs}
\usepackage{graphics, epsfig, here, graphicx, setspace}
\usepackage{bbm, color}
\usepackage[all]{xy}
\definecolor{darkgreen}{rgb}{0,0.5,0}
\definecolor{darkred}{rgb}{0.5,0,0}
\definecolor{darkblue}{rgb}{0,0,0.5}
\definecolor{rb1}{rgb}{1,0,0}
\definecolor{rb2}{rgb}{.75,0,.25}
\definecolor{rb3}{rgb}{.5,0,.5}
\definecolor{rb4}{rgb}{.25,0,.75}
\definecolor{rb5}{rgb}{0,0,1}

\newcommand{\boX}{\hfill $\Box$}
\newcommand{\bN}{\mathbb{N}}

\newcommand{\bR}{\mathbb{R}}

\newcommand{\cA}{\mathcal{A}}
\newcommand{\cP}{\mathcal{P}}

\newcommand{\cJ}{\mathcal{J}}

\newcommand{\cD}{\mathcal{D}}

\newcommand{\cO}{\mathcal{O}}

\newcommand{\eX}{\exists\,}
\newcommand{\fA}{\,\forall\,}
\newcommand{\ol}[1]{{\overline{#1}}}
\newcommand{\dS}{\displaystyle}
\newcommand{\ssS}{\scriptscriptstyle}
\newcommand{\uNM}{{\ssS \! \! N\! M}}

\newcommand{\bc}[1]{{\left\{#1\right\}}}
\newcommand{\bq}[1]{{\left[#1\right]}}
\newcommand{\bp}[1]{{\left(#1\right)}}
\newcommand{\bl}[1]{{\left|#1\right|}}
\newcommand{\bt}[1]{{\left\langle#1\right\rangle}}

\newcommand{\norms}[1]{{\|#1\|}}
\newcommand{\set}[2]{{\bc{#1\,\,:\,\,#2}}}

\newcommand{\teQ}{\triangleq}

\newcommand{\rM}[1]{{\mathrm{#1}}}

\newcommand{\tb}[1]{{\textbf{#1}}}

\newcommand{\ep}{\epsilon}
\newcommand{\vphi}{\varphi}
\newcommand{\Dx}{\dot{x}}

\newcommand{\DPhi}{\dot{\Phi}}

\renewcommand{\O}{\Omega}

\newcommand{\w}{\omega}
\newcommand{\cF}{\mathcal{F}}

\newcommand{\st}{\,\,\mathrm{ s.t. }\,\,}

\newcommand{\anD}{\quad\text{ and }\quad}
\newcommand{\eP}{\varepsilon}
\newcommand{\evalat}[3]{{\Bigl.#1 \Bigr|_{#2}^{#3}}}
\newcommand{\evalatb}[3]{{\bigl.#1 \bigr|_{#2}^{#3}}}

\newcommand{\tw}[1]{{\widetilde{#1}}}

\newcommand{\half}{\frac{1}{2}}

\newcommand{\rd}{\mathrm{d}}
\newcommand{\ddx}[1]{{\frac{\rd }{\rd #1}}}
\newcommand{\ppx}[1]{{\frac{\partial }{\partial #1}}}

\newcommand{\dxdy}[2]{{\frac{\rd #1}{\rd #2}}}

%% file: nm_ent_s0.bbl
\begin{thebibliography}{10}

\bibitem{blekhman88}
I.~Blekhman.
\newblock {\em Synchronization in science and technology}.
\newblock ASME Press translations, New York, 1988.

\bibitem{pikovsky01}
A.~Pikovsky, M.~Rosenblum, and J.~Kurths.
\newblock {\em Synchronization: {A} {U}niversal {C}oncept in {N}onlinear
  {S}cience}.
\newblock Cambridge University Press, 2001.

\bibitem{strogatz01}
S.~Strogatz.
\newblock {\em Nonlinear {D}ynamics And {C}haos: {W}ith {A}pplications To
  {P}hysics, {B}iology, {C}hemistry, And {E}ngineering}.
\newblock Studies in nonlinearity. Westview Press, 1 edition, 2001.

\bibitem{granada09enz}
A.~Granada, R.~M. Hennig, B.~Ronacher, A.~Kramer, and H.~Herzel.
\newblock Phase response curves: elucidating the dynamics of coupled
  oscillators.
\newblock {\em Methods in Enzymology}, 454:1--27, 2009.

\bibitem{izhikevich07}
E.~Izhikevich.
\newblock {\em Dynamical {S}ystems in {N}euroscience}.
\newblock Neuroscience. MIT Press, 2007.

\bibitem{mrosovsky80}
N.~Mrosovsky.
\newblock Circannual cycles in golden-mantled ground squirrels: phase shift
  produced by low temperatures.
\newblock {\em Journal of Comparative Physiology A: Neuroethology, Sensory,
  Neural, and Behavioral Physiology}, 136(4):349--353, 1980.

\bibitem{doyle06}
F.~J. Doyle, R.~Gunawan, N.~Bagheri, H.~Mirsky, and T.~L. To.
\newblock Circadian rhythm: A natural, robust, multi-scale control system.
\newblock {\em Computers \& chemical engineering}, 30(10):1700--1711, 2006.

\bibitem{mcclung11}
{McClung, C. R. et. al}.
\newblock The genetics of plant clocks.
\newblock {\em Advances in genetics}, 74:105--139, 2011.

\bibitem{gonze05}
D.~Gonze, S.~Bernard, C.~Waltermann, A.~Kramer, and H.~Herzel.
\newblock Spontaneous synchronization of coupled circadian oscillators.
\newblock {\em Biophysical Journal}, 89(1):120--129, 2005.

\bibitem{coombes05}
S.~Coombes and P.C. Bressloff.
\newblock {\em Bursting: The genesis of rhythm in the nervous system}.
\newblock World Scientific Publishing Company Incorporated, 2005.

\bibitem{hartbauer05}
M.~Hartbauer, S.~Kratzer, K.~Steiner, and H.~R{\"o}mer.
\newblock Mechanisms for synchrony and alternation in song interactions of the
  bushcricket mecopoda elongata (tettigoniidae: Orthoptera).
\newblock {\em Journal of Comparative Physiology A: Neuroethology, Sensory,
  Neural, and Behavioral Physiology}, 191(2):175--188, 2005.

\bibitem{strogatz00}
S.~H. Strogatz.
\newblock From {K}uramoto to {C}rawford: exploring the onset of synchronization
  in populations of coupled oscillators.
\newblock {\em Physica D}, 143(1-4):1--20, 2000.

\bibitem{hanson78}
F.~Hanson.
\newblock Comparative studies of firefly pacemakers.
\newblock {\em Federation Proceedings}, 38(8):2158--2164, 1978.

\bibitem{ermentrout84}
G.~Ermentrout and J.~Rinzel.
\newblock Beyond a pacemaker's entrainment limit: phase walk-through.
\newblock {\em American Journal of Physiology - Regulatory, Integrative and
  Comparative Physiology}, 246(1), 1984.

\bibitem{glass88}
L.~Glass and M.C. Mackey.
\newblock {\em From clocks to chaos: The rhythms of life}.
\newblock Princeton University Press, 1988.

\bibitem{demir97}
S.~Demir, R.~Butera, A.~DeFranceschi, J.~Clark, and J.~Byrne.
\newblock Phase {S}ensitivity and {E}ntrainment in a {M}odeled {B}ursting
  {N}euron.
\newblock {\em Biophysical Journal}, 72:579--594, 1997.

\bibitem{berke04}
J.~D. Berke, M.~Okatan, J.~Skurski, and H.~B. Eichenbaum.
\newblock Oscillatory {E}ntrainment of {S}triatal {N}eurons in {F}reely
  {M}oving {R}ats.
\newblock {\em Neuron}, 43:883--896, 2004.

\bibitem{sirota08}
A.~Sirota, S.~Montgomery, S.~Fujisawa, Y.~Isomura, M.~Zugaro, and G.~Buzs\'aki.
\newblock Entrainment of neocortical neurons and gamma oscillations by the
  hippocampal theta rhythm.
\newblock {\em Neuron}, 60:683--697, 2008.

\bibitem{kuramoto84}
Y.~Kuramoto.
\newblock {\em Chemical {O}scillations, {W}aves, and {T}urbulence}.
\newblock Springer, New York, 1984.

\bibitem{aronson86}
D.~Aronson, R.~McGehee, I.~Kevrekidis, and R.~Aris.
\newblock Entrainment regions for periodically forced oscillators.
\newblock {\em Physical Review A}, 33(3):2190--2192, 1986.

\bibitem{lev89}
O.~Lev, A.~Wolfberg, L.~M. Pismen, and M.~Sheintuch.
\newblock The structure of complex behavior in anodic nickel dissolution.
\newblock {\em Journal of Phys. Chem.}, 93:1661--1666, 1989.

\bibitem{granada09fast}
A.~E. Granada and H.~Herzel.
\newblock How to achieve fast entrainment? the timescale to synchronization.
\newblock {\em PLoS One}, 4(9):e7057, 2009.

\bibitem{harada10}
T.~Harada, H.~Tanaka, M.~Hankins, and I.~Kiss.
\newblock Optimal waveform for the entrainment of a weakly forced oscillator.
\newblock {\em Physical Review Letters}, 105(8), 2010.

\bibitem{vosko10}
A.~M. Vosko, C.~S. Colwell, and A.~Y. Avidan.
\newblock Jet lag syndrome: circadian organization, pathophysiology, and
  management strategies.
\newblock {\em Nature}, 2:187--198, 2010.

\bibitem{sthilaire12}
M.~A. St.~Hilaire, J.~J. Gooley, S.~B.~S. Khalsa, R.~E. Kronauer, C.~A.
  Czeisler, and S.~W. Lockley.
\newblock Human phase response curve to a 1 h pulse of bright white light.
\newblock {\em The Journal of Physiology}, 2012.

\bibitem{kiss08}
I.~Z. Kiss, M.~Quigg, S.~H.~C. Chun, H.~Kori, and J.~L. Hudson.
\newblock Characterization of synchronization in interacting groups of
  oscillators: application to seizures.
\newblock {\em Biophysical Journal}, 94(3):1121--1130, 2008.

\bibitem{good09}
L.~Good.
\newblock Control of synchronization of brain dynamics leads to control of
  epileptic seizures in rodents.
\newblock {\em International Journal of Neural Systems}, 19(3):173--196, 2009.

\bibitem{hofmann11}
L.~Hofmann, M.~Ebert, P.A. Tass, and C.~Hauptmann.
\newblock Modified pulse shapes for effective neural stimulation.
\newblock {\em Frontiers in Neuroengineering}, 4, 2011.

\bibitem{strauss05}
D.~J. Strauss, W.~Delb, R.~D'Amelio, and P.~Falkai.
\newblock Neural synchronization stability in the tinnitus decompensation.
\newblock In {\em Neural Engineering, 2005. Conference Proceedings. 2nd
  International IEEE EMBS Conference on}, pages 186--189. IEEE, 2005.

\bibitem{naqvi12}
T.~Z. Naqvi and D.~C. Winter.
\newblock Optimization of pacemaker settings, February~7 2012.
\newblock {US} Patent 8,112,150.

\bibitem{zalalutdinov03}
M.~Zalalutdinov, K.~Aubin, A.~Zehnder, R.~Hand, H.~Craighead, J.~Parpia, and
  B.~Houston.
\newblock Frequency entrainment for micromechanical oscillator.
\newblock {\em Applied Physics Letters}, 83(16):3281--3283, 2003.

\bibitem{feng08}
X.~L. Feng, C.~J. White, A.~Hajimiri, and M.~L. Roukes.
\newblock A self-sustaining ultrahigh-frequency nanoelectromechanical
  oscillator.
\newblock {\em Nature Nanotechnology}, 3(6):342--346, 2008.

\bibitem{barnes11}
A.~C. Barnes, R.~C. Roberts, N.~C. Tien, C.~A. Zorman, and P.~X.~L. Feng.
\newblock Silicon carbide (sic) membrane nanomechanical resonators with
  multiple vibrational modes.
\newblock In {\em Solid-State Sensors, Actuators and Microsystems Conference
  (TRANSDUCERS), 2011 16th International}, pages 2614--2617. IEEE, 2011.

\bibitem{hoppensteadt00}
F.~C. Hoppensteadt and E.~M. Izhikevich.
\newblock Synchronization of laser oscillators, associative memory, and optical
  neurocomputing.
\newblock {\em Physical Review E}, 62(3):4010--4013, 2000.

\bibitem{fischer00}
I.~Fischer, Y.~Liu, and P.~Davis.
\newblock Synchronization of chaotic semiconductor laser dynamics on
  subnanosecond time scales and its potential for chaos communication.
\newblock {\em Physical Review A}, 62, 2000.

\bibitem{wacker11}
M.~Wacker and H.~Witte.
\newblock On the stability of the n: m phase synchronization index.
\newblock {\em Biomedical Engineering, IEEE Transactions on}, 58(2):332--338,
  2011.

\bibitem{chen09}
L.~Chen, C.~Qiu, and HB~Huang.
\newblock Synchronization with on-off coupling: Role of time scales in network
  dynamics.
\newblock {\em Physical Review E}, 79(4):045101, 2009.

\bibitem{chavez05}
M.~Chavez, C.~Adam, V.~Navarro, S.~Boccaletti, and J.~Martinerie.
\newblock On the intrinsic time scales involved in synchronization: a
  data-driven approach.
\newblock {\em Chaos: An Interdisciplinary Journal of Nonlinear Science},
  15(2):023904--023904, 2005.

\bibitem{zlotnik13prl}
A.~Zlotnik, Y.~Chen, I.~Z. Kiss, H.-A. Tanaka, and J.-S. Li.
\newblock Optimal waveform for fast entrainment of weakly forced nonlinear
  oscillators.
\newblock {\em Phys. Rev. Lett.}, 111:024102, Jul 2013.

\bibitem{ermentrout96}
B.~Ermentrout.
\newblock Type {I} {M}embranes, {P}hase {R}esetting {C}urves, and {S}ynchrony.
\newblock {\em Neural Computation}, 8(5):979--1001, 1996.

\bibitem{galan05}
R.~F. Gal{\'a}n, G.~B. Ermentrout, and N.~N. Urban.
\newblock Efficient estimation of phase-resetting curves in real neurons and
  its significance for neural-network modeling.
\newblock {\em Physical Review Letters}, 94(15):158101, 2005.

\bibitem{tokuda07}
I.~T. Tokuda, S.~Jain, I.~Z. Kiss, and J.~L. Hudson.
\newblock Inferring phase equations from multivariate time series.
\newblock {\em Physical Review Letters}, 99(6):64101, 2007.

\bibitem{hoppensteadt99}
F.~Hoppensteadt and E.~Izhikevich.
\newblock Oscillatory {N}eurocomputers with {D}ynamic {C}onnectivity.
\newblock {\em Physical Review Letters}, 82(14), 1999.

\bibitem{kiss02}
I.~Z. Kiss, Y.~Zhai, and J.~Hudson.
\newblock Emerging coherence in a population of chemical oscillators.
\newblock {\em Science}, 296:1676--1678, 2002.

\bibitem{nakata09}
S.~Nakata, K.~Miyazaki, S.~Izuhara, H.~Yamaoka, and D.~Tanaka.
\newblock Arnold {T}ongue of {E}lectrochemical {N}onlinear {O}scillators.
\newblock {\em Journal of Physical Chemistry A}, 113:6876--6879, 2009.

\bibitem{moehlis06}
H.~Moehlis, E.~Brown, and H.~Rabitz.
\newblock Optimal inputs for phase models of spiking neurons.
\newblock {\em Journal of Computational and Nonlinear Dynamics}, 1:358--367,
  2006.

\bibitem{dasanayake11}
I.~Dasanayake and J.-S. Li.
\newblock Optimal design of minimum-power stimuli for phase models of neuron
  oscillators.
\newblock {\em Phys. Rev. E}, 83:061916, 2011.

\bibitem{dasanayake11cdc}
I.~Dasanayake and J.-S. Li.
\newblock Constrained minimum-power control of spiking neuron oscillators.
\newblock In {\em IEEE CDC}, pages 3694--3699. IEEE, 2011.

\bibitem{danzl10}
P.~Danzl, A.~Nabi, and J.~Moehlis.
\newblock Charge-balanced spike timing control for phase models of spiking
  neurons.
\newblock {\em Discrete and Continuous Dynamical Systems}, 28(4):1413--1435,
  2010.

\bibitem{dasanayake11cb}
I.~Dasanayake and Li~J.-S.
\newblock Charge-balanced minimum-power controls for spiking neuron
  oscillators.
\newblock {\em arXiv:1109.3798}, 2011.

\bibitem{nabi11}
A.~Nabi and J.~Moehlis.
\newblock Single input optimal control for globally coupled neuron networks.
\newblock {\em J. Neural Eng.}, 8:065008, 2011.

\bibitem{schaus06}
M.~J. Schaus and J.~Moehlis.
\newblock On the response of neurons to sinusoidal current stimuli: Phase
  response curves and phase-locking.
\newblock In {\em Proceedings of the 45th IEEE Conference on Decision \&
  Control}, pages 2376--2381, San Diego, CA, December 2006.

\bibitem{zlotnik11}
A.~Zlotnik and J.~Li.
\newblock Optimal asymptotic entrainment of phase-reduced oscillators.
\newblock In {\em 2011 ASME Dynamic Systems and Control Conference}, volume~1,
  pages 479--484, Arlington, VA, October 2011.

\bibitem{zlotnik12jne}
A.~Zlotnik and J.-S. Li.
\newblock Optimal entrainment of neural oscillator ensembles.
\newblock {\em J. Neural Eng.}, 9(4):046015, 2012.

\bibitem{cunningham51}
W.~J. Cunningham.
\newblock The growth of subharmonic oscillations.
\newblock {\em The Journal of the Acoustical Society of America}, 23:418, 1951.

\bibitem{ermentrout81}
G.~B. Ermentrout.
\newblock {n:m} phase-locking of weakly coupled oscillators.
\newblock {\em Journal of Mathematical Biology}, 12:327--342, 1981.

\bibitem{guevara82}
M.~Guevara and L.~Glass.
\newblock Phase {L}ocking, {P}eriod {D}oubling {B}ifurcations and {C}haos in a
  {M}athematical {M}odel of a {P}eriodically {D}riven {O}scillator.
\newblock {\em Journal of Mathematical Biology}, 14:1--23, 1982.

\bibitem{daryoush89}
A.~S. Daryoush, T.~Berceli, R.~Saedi, P.~R. Herczfeld, and A.~Rosen.
\newblock Theory of subharmonic synchronization of nonlinear oscillators.
\newblock In {\em Microwave Symposium Digest, 1989., IEEE MTT-S International},
  pages 735--738. IEEE, 1989.

\bibitem{storti88}
D.~Storti and R.H. Rand.
\newblock Subharmonic entrainment of a forced relaxation oscillator.
\newblock {\em International journal of non-linear mechanics}, 23(3):231--239,
  1988.

\bibitem{tass98}
P.~Tass, MG~Rosenblum, J.~Weule, J.~Kurths, A.~Pikovsky, J.~Volkmann,
  A.~Schnitzler, and H.J. Freund.
\newblock Detection of n: m phase locking from noisy data: application to
  magnetoencephalography.
\newblock {\em Physical Review Letters}, 81(15):3291--3294, 1998.

\bibitem{honey07}
C.~J. Honey, R.~K{\"o}tter, M.~Breakspear, and O.~Sporns.
\newblock Network structure of cerebral cortex shapes functional connectivity
  on multiple time scales.
\newblock {\em Proceedings of the National Academy of Sciences},
  104(24):10240--10245, 2007.

\bibitem{hunter03}
J.~Hunter and J.~Milton.
\newblock Amplitude and {F}requency {D}ependence of {S}pike {T}iming:
  {I}mplications for {D}ynamic {R}egulation.
\newblock {\em Journal of Neurophysiology}, 90:387--394, 2003.

\bibitem{kriellaars94}
D.~J. Kriellaars, R.~M. Brownstone, B.~R. Noga, and L.~M. Jordan.
\newblock Mechanical entrainment of fictive locomotion in the decerebrate cat.
\newblock {\em Journal of neurophysiology}, 71(6):2074--2086, 1994.

\bibitem{large96}
E.~W. Large.
\newblock Modeling beat perception with a nonlinear oscillator.
\newblock In {\em Proceedings of the 18th Annual Conference of the Cognitive
  Science Society}, page 420, 1996.

\bibitem{clayton05}
M.~Clayton, R.~Sager, and U.~Will.
\newblock In time with the music: The concept of entrainment and its
  significance for ethnomusicology.
\newblock In {\em European meetings in ethnomusicology}, volume~11, pages
  3--142, 2005.

\bibitem{nozaradan11}
S.~Nozaradan, I.~Peretz, M.~Missal, and A.~Mouraux.
\newblock Tagging the neuronal entrainment to beat and meter.
\newblock {\em The Journal of Neuroscience}, 31(28):10234--10240, 2011.

\bibitem{tass12}
P.~A. Tass, I.~Adamchic, H.~J. Freund, T.~von Stackelberg, and C.~Hauptmann.
\newblock Counteracting tinnitus by acoustic coordinated reset neuromodulation.
\newblock {\em Restorative neurology and neuroscience}, 30(2):137--159, 2012.

\bibitem{lajoie11}
G.~Lajoie and E.~Shea-Brown.
\newblock Shared inputs, entrainment, and desynchrony in elliptic bursters:
  from slow passage to discontinuous circle maps.
\newblock {\em SIAM Journal on Applied Dynamical Systems}, 10(4):1232--1271,
  2011.

\bibitem{gutierrez11}
R.~Guti{\'e}rrez, A.~Amann, S.~Assenza, J.~G{\'o}mez-Gardenes, V.~Latora, and
  S.~Boccaletti.
\newblock Emerging meso-and macroscales from synchronization of adaptive
  networks.
\newblock {\em Physical Review Letters}, 107(23):234103, 2011.

\bibitem{schafer98}
C.~Sch{\"a}fer, M.G. Rosenblum, J.~Kurths, and H.H. Abel.
\newblock Heartbeat synchronized with ventilation.
\newblock {\em Nature}, 392:239--240, 1998.

\bibitem{cajochen13}
C.~Cajochen, S.~Altanay-Ekici, M.~M{\"u}nch, S.~Frey, V.~Knoblauch, and
  A.~Wirz-Justice.
\newblock Evidence that the lunar cycle influences human sleep.
\newblock {\em Current Biology}, 2013.

\bibitem{foster08}
R.~G. Foster and T.~Roenneberg.
\newblock Human responses to the geophysical daily, annual and lunar cycles.
\newblock {\em Current biology}, 18(17):R784--R794, 2008.

\bibitem{maffezzoni10}
P.~Maffezzoni, D.~D'Amore, S.~Daneshgar, and MP~Kennedy.
\newblock Estimating the locking range of analog dividers through a
  phase-domain macromodel.
\newblock In {\em Circuits and Systems (ISCAS), Proceedings of 2010 IEEE
  International Symposium on}, pages 1535--1538. IEEE, 2010.

\bibitem{takano07}
K.~Takano, M.~Motoyoshi, and M.~Fujishima.
\newblock 4.8 ghz cmos frequency multiplier with subharmonic pulse-injection
  locking.
\newblock In {\em Solid-State Circuits Conference, 2007. ASSCC'07. IEEE Asian},
  pages 336--339. IEEE, 2007.

\bibitem{zarroug95}
A.~Zarroug, PS~Hall, and M.~Cryan.
\newblock Active antenna phase control using subharmonic locking.
\newblock {\em Electronics Letters}, 31(11):842--843, 1995.

\bibitem{bolanos05}
F.~Bola\~{n}os.
\newblock Measurement and analysis of subharmonics and other distortions in
  compression drivers.
\newblock volume~1, 2005.

\bibitem{noris12}
J.~W. Noris.
\newblock Nonlinear dynamical behavior of a moving voice coil.
\newblock In {\em 105th AES Convention, San Francisco}, volume~1, 1988.

\bibitem{hodgkin52}
A.~Hodgkin and A.~Huxley.
\newblock A quantitative description of membrane current and its application to
  conduction and excitation in nerve.
\newblock {\em The Journal of Physiology}, 117(4), 1952.

\bibitem{kiss05}
I.~Z. Kiss, Y.~M. Zhai, and J.~L. Hudson.
\newblock Predicting mutual entrainment of oscillators with experiment-based
  phase models.
\newblock {\em Phys. Rev. Lett.}, 94(24):248301, 2005.

\bibitem{brown04}
E.~Brown, J.~Moehlis, and P.~Holmes.
\newblock On the {P}hase {R}eduction and {R}esponse {D}ynamics of {N}eural
  {O}scillator {P}opulations.
\newblock {\em Neural Computation}, 16(4):673--715, 2004.

\bibitem{efimov10}
D.~Efimov and T.~Raissi.
\newblock Phase resetting control based on direct phase response curve.
\newblock In {\em Preprints of the 8th IFAC Symposium on Nonlinear Control
  Systems}, pages 332--337, Bologna, September 2010.

\bibitem{efimov11}
D.~Efimov.
\newblock Phase resetting control based on direct phase response curve.
\newblock {\em Journal of mathematical biology}, 63(5):855--879, 2011.

\bibitem{efimov09}
D.~Efimov, P.~Sacr\'e, and R.~Sepulchre.
\newblock Controlling the {P}hase of an {O}scillator: {A} {P}hase {R}esponse
  {C}urve {A}pproach.
\newblock In {\em Joint 48th Conference on Decision and Control}, pages
  7692--7697, December 2009.

\bibitem{perko90}
L.~Perko.
\newblock {\em Differential equations and dynamical systems}.
\newblock Texts in applied mathematics. Springer, 2 edition, 1990.

\bibitem{kelley04}
W.~Kelley and A.~Peterson.
\newblock {\em The {T}heory of {D}ifferential {E}quations, {C}lassical and
  {Q}ualitative}.
\newblock Pearson, 2004.

\bibitem{aprille72}
T.~Aprille and T.~Trick.
\newblock A computer algorithm to determine the steady-state response of
  nonlinear oscillators.
\newblock {\em IEEE Transactions on Circuit Theory}, 19(4):354--360, 1972.

\bibitem{khalil02}
H.~Khalil.
\newblock {\em Nonlinear {S}ystems}.
\newblock Prentice Hall, 3 edition, 2002.

\bibitem{peressini00}
A.~Peressini, F.~Sullivan, and J.~Uhl.
\newblock {\em Mathematics of {N}onlinear {P}rogramming}.
\newblock Springer, 2000.

\bibitem{govaerts06}
W.~Govaerts and B.~Sautois.
\newblock Computation of the {P}hase {R}esponse {C}urve: {A} {D}irect
  {N}umerical {A}pproach.
\newblock {\em Neural Computation}, 18(4):817--847, 2006.

\bibitem{ermentrout02}
B.~Ermentrout.
\newblock {\em Simulating, {A}nalyzing, and {A}nimating {D}ynamical {S}ystems:
  A {G}uide to {XPPAUT} for {R}esearchers and {S}tudents}.
\newblock SIAM, 2002.

\bibitem{malkin49}
I.~Malkin.
\newblock {\em Methods of {P}oincare and {L}iapunov in the theory of nonlinear
  oscillations}.
\newblock Gostexizdat, Moscow, 1949.

\bibitem{kornfeld82}
I.~Kornfeld, S.~Fomin, and Y.~Sinai.
\newblock {\em Ergodic theory: Differentiable Dynamical Systems}, volume 245 of
  {\em Grund. Math. Wissens.}
\newblock Springer-Verlag, 1982.

\bibitem{hoppensteadt97}
F.~Hoppensteadt and E.~Izhikevich.
\newblock {\em Weakly connected neural networks}.
\newblock Springer-Verlag, New Jersey, 1997.

\bibitem{arnold61}
V.~I. Arnol'd.
\newblock Small denominators. i. mapping the circle onto itself.
\newblock {\em Izvestiya Rossiiskoi Akademii Nauk. Seriya Matematicheskaya},
  25(1):21--86, 1961.

\bibitem{gelfand00}
I.~Gelfand and S.~Fomin.
\newblock {\em Calculus of {V}ariations}.
\newblock Dover, 2000.

\bibitem{kiss07}
I.~Z. Kiss, C.~G. Rusin, H.~Kori, and J.~L. Hudson.
\newblock Engineering complex dynamical structures: sequential patterns and
  desynchronization.
\newblock {\em Science}, 316(5833):1886--1889, 2007.

\bibitem{uhlhaas06}
P.~J. Uhlhaas and W.~Singer.
\newblock Neural synchrony in brain disorders: relevance for cognitive
  dysfunctions and pathophysiology.
\newblock {\em Neuron}, 52(1):155--168, 2006.

\bibitem{gooley10}
J.~J. Gooley, S.~M. Rajaratnam, G.~C. Brainard, R.~E. Kronauer, C.~A. Czeisler,
  and S.~W. Lockley.
\newblock Spectral responses of the human circadian system depend on the
  irradiance and duration of exposure to light.
\newblock {\em Science Translational Medicine}, 2(31):31ra33, 2010.

\end{thebibliography}
